\documentclass[11pt]{article}
\usepackage{amssymb,amsthm,amsfonts,hhline,color}
\usepackage{longtable,ltcaption}

\def\F{{\bf F}}
\def\Z{{\bf Z}}

\def\Q{{\bf Q}}

\def\K32{K3^{[2]}}
\def\K3n{K3^{[n]}}
\def\Co{{\rm Co}_0}
\def\mon{2^{12}{:}M_{24}}
\def\M24{M_{24}}
\def\M23{M_{23}}
\def\mcl{{\rm McL}}
\def\S{{\mathcal{S}}}

\def\rk{{\rm rk}\,}

\def\pf{\noindent{\bf Proof:\ }}
\def\qed{\hfill\framebox[2.5mm][t1]{\phantom{x}}}

\title{The $290$ fixed-point sublattices of the Leech lattice}

\author{Gerald H\"ohn\\
Department of Mathematics, Kansas State University
\\  Geoffrey Mason\thanks{Supported by the NSF}\\
Department of Mathematics, University of California at Santa Cruz}

\date{May, 2015}

\begin{document}

\bibliographystyle{amsalpha}

\theoremstyle{plain}
\newtheorem{thm}{Theorem}[section]
\newtheorem{prop}[thm]{Proposition}
\newtheorem{lem}[thm]{Lemma}
\newtheorem{cor}[thm]{Corollary}
\newtheorem{rem}[thm]{Remark}
\newtheorem{conj}[thm]{Conjecture}

\newtheorem{introthm}{Theorem}
\renewcommand\theintrothm{\Alph{introthm}}

\theoremstyle{definition}
\newtheorem{defi}[thm]{Definition}

\renewcommand{\baselinestretch}{1.2}

\def\pf{\noindent{\bf Proof:\ }}
\def\qed{\hfill\framebox[2.5mm][t1]{\phantom{x}}}

\maketitle

\begin{abstract}
\noindent
We determine the orbits of fixed-point sublattices of the Leech lattice with respect to the action
of the Conway group $\Co$.\ There are 290 such orbits.\  Detailed information about
these lattices, the corresponding coinvariant lattices, and the stabilizing subgroups, is tabulated in several tables.

\end{abstract}

%%%%%%%%%%%%%%%%%%%%%%%%%%%%%%%%%%%%%%%%%%%%%%%%%%%%%%%%%%%%%%%%%%%%%%%%%%%%%%%%%%%%%%%%%%%%%%%%%%%%%%%%%%%%%%%%%

\section{Introduction}

The {\it Leech lattice} $\Lambda$ is the unique positive-definite, even, 
unimodular lattice of rank~$24$ without roots~\cite{Leech, Conway2}.\
It may also be characterized as the most densely packed lattice in 
dimension~$24$~\cite{CoKu}.\
The group of isometries of $\Lambda$ is the {\it Conway group\/} $\Co$~\cite{Conway1}.\
For a subgroup $H\subseteq \Co$ we set
\begin{eqnarray*}
\Lambda^H & =  &\{ v\in\Lambda\ \mid \ hv=v\ \mbox{\ for all\ } h\in H\}.
\end{eqnarray*}
We call such a sublattice of $\Lambda$ a \emph{fixed-point sublattice.\/}\
Let $\mathcal{F}$ be the set of all fixed-point sublattices of $\Lambda$.\ 
The Conway group acts by translation on $\mathcal{F}$, because if $g\in \Co$,
then $g \Lambda^G  = \Lambda^{gHg^{-1}}.$\ In this note, we classify the $\Co$-orbits of fixed-point sublattices.\
We will prove:
\begin{thm}\label{mainthm}
Under the action of $\Co$, there are exactly $290$ orbits on the set of fixed-point sublattices of $\Lambda$.
\end{thm}

The purpose of the present note is not merely to enumerate the orbits of fixed-point sublattices, 
but to provide in addition a detailed analysis of their properties.\
In particular, this includes the \emph{stabilizers\/} $G$, which are the (largest) subgroups of $\Co$ 
that stabilize a given fixed-point sublattice pointwise.\
Information about the orbits of fixed-point lattices and their fixing groups is given in 
Table~\ref{maintable} in Section~\ref{SSA}.  Based on the theory that we present in Sections~\ref{basicstuff}
and~\ref{construction}, this information was obtained  by relying on extensive computer calculations using
the computer algebra system MAGMA~\cite{magma}.\ We shall say more about this in due course.

\medskip

There are a number of reasons that make the classification of fixed-point lattices desirable.\  
The group quotient ${\rm Co}_1=\Co/\{\pm 1\}$  is one of the $26$ sporadic simple groups.\ It contains
$11$ additional sporadic groups, $9$ of which can be described in terms of lattice stabilizers.\ 
Although these particular realizations have been known for a long time, the complete picture that we provide is new.

The Leech lattice is also the starting point of the construction of interesting
{\it vertex operator algebras\/}~\cite{Bo-ur,FLM}
and {\it generalized Kac-Moody Lie algebras.\/}\ Such Kac-Moody Lie algebras have root 
lattices that can often be described in terms of fixed-point lattices inside $\Lambda$~\cite{Scheit}, 
and the associated denominator identities provide Moonshine for the corresponding subgroups~\cite{Bo-Moon}.\

\medskip
The {\it geometry of $K3$ surfaces\/} and certain {\it hyperk\"ahler manifolds\/} $X$,
over both the field of complex numbers and in finite characteristic,
is controlled (using Torelli-type theorems) by lattices related to $\Lambda$.\ 
In this way, symmetry groups of $X$ can be mapped into $\Co$, and properties of the fixed-point lattices 
control which groups may appear.\
See~\cite{Nikulin,Mukai,Kondo,DoKe} for $K3$ surfaces and \cite{Mon-thesis,Huy-conway,K32} for other 
hyperk\"ahler manifolds.\ 

\smallskip

Much of the impetus for studying the finite symmetry groups of such manifolds, and recent developments in the related area of
\emph{Mathieu Moonshine\/}~\cite{EOT}, came from the well-known theorem of Mukai~\cite{Mukai}.\ 
This states that a finite group $G$ of symplectic automorphisms of a K3 surface is isomorphic to a subgroup of the Mathieu group $M_{23}$ 
with at least five orbits in its natural permutation representation on $24$ letters; furthermore, there are just $11$ subgroups 
(up to isomorphism) which are maximal subject to these conditions.\ A typical application of our results 
leads to a simplified approach to this theorem.\ Indeed, lattice-theoretic arguments~\cite{Kondo,Conway3} show that $G$ 
can be embedded into $\Co$ in such a way that ${\rm rk}\,\Lambda^G\geq 5$  
and $\alpha(\Lambda^G):={\rm rk}\,  \Lambda^G - {\rm rk}\,A_{\Lambda^G} \geq 2$ (see Section~\ref{basicstuff} for notation).\
The containment $G\subseteq M_{23}$ follows immediately from Table~\ref{maintable}, moreover the $11$ maximal such groups are 
those $G$ in Table~\ref{maintable} with ${\rm rk}\,  \Lambda^G=5$ and $\alpha(\Lambda^G)\ge 2$.
The advantage of this approach compared to that of Kondo (loc.\ cit), who initiated the lattice-theoretic approach, 
is that a case-by-case analysis of the $23$ Niemeier lattices with roots is reduced to an analysis of the Leech lattice alone. 

The extension of Mukai's theorem to more general contexts is currently an active research area,
and it is widely expected that knowledge of the stabilizers $G$ with $\Lambda^G \geq 4$ will eventually lead to the classification of all
finite symplectic automorphism groups of hyperk\"ahler manifolds of type $K3^{[n]}$ (cf.~\cite{Mo2,K32}) and symplectic autoequivalences
of derived categories of sheaves on K3 surfaces~\cite{AM-K3,GHV,Huy-conway}.

\smallskip

Another application of our tables is to the study of symmetries of the extremal vertex operator superalgebra $V_{D_{12}^+}$ of central
charge~$12$\/~\cite{Ho-thesis}.\ This SVOA has a unique $N\!\!=\!\!1$ super conformal structure stabilized by the Conway group~\cite{Duncan-thesis}.\ 
By considering subgroups $G\subseteq \Co$ fixing subspaces of $\Lambda\otimes {\bf R}$ of various
dimensions~\cite{Kachru1,Kachru2,Duncan-Jacobi} one can define Jacobi forms of weight~$0$ 
(and some index), leading to various moonshine phenomena for these $G$.\ The groups $G$ which can occur can directly be read off 
from Table~\ref{maintable}.\ Only a few examples of such groups have been known and studied before.\ The resulting moonshine functions 
are also closely related to the equivariant elliptic genera of hyperk\"ahler manifolds of type $K3^{[n]}$ (cf.~\cite{CH,K32}) and to
to Mathieu Moonshine and its generalization called Umbral Moonshine~\cite{Umbral}; see~\cite{Duncan-Jacobi}.

\medskip

Notable past work on sublattices of $\Lambda$ includes that of 
Curtis~\cite{Curtis}, who introduced the very useful class of \emph{${\S}$-lattices} and classified their orbits.\ 
We make use of these ideas in the present paper (cf.\ Section~\ref{construction}).\
Harada and Lang~\cite{HaLa} considered the orbits of fixed-point lattices for \emph{cyclic\/} subgroups of $\Co$.\
Fixed-point lattices related to $K3$ surfaces were classified by Hashimoto~\cite{Ha}, and the present authors handled those 
for hyperk\"ahler manifolds of type $K3^{[2]}$~\cite{K32}.\
Additional information can also be found in the Atlas of finite groups~\cite{Atlas}.

\medskip

The present work is based on lattice-theoretic arguments and group-theoretical computations in $\Co$. \

\smallskip

The Conway group $\Co$ is presently too large to permit computation of its complete subgroup lattice.\
(Such a calculation would allow us to list all orbits of fixed-point lattices directly.)\ 
Even for the monomial subgroup $\mon\subseteq \Co$,
the number of conjugacy classes is huge (of order $10^7$--$10^8$) and we have thus far
been unable to determine them all.\ For the purposes of the present work, however,  it
is enough to know the conjugacy classes of non-$2$-groups inside $\mon$ and these were already computed in~\cite{K32}.\  
(There are $279,343$ of them.)\

\bigskip

The paper is organized as follows.\ Section~\ref{basicstuff} summarizes some general properties
of group actions on lattices.\ In Section~\ref{construction}, we describe our method to
determine the $290$ orbits, while Section~\ref{SSA} contains detailed information about the $290$ fixed-point lattices.
We also discuss several interesting properties of some of the resulting lattices.

\medskip

Finally, we mention that the corresponding problem of classification of fixed-point lattices and 
stabilizer subgroups in the case of the $E_8$-root lattice and its attendant Weyl group is also of interest.\
For the convenience of the interested reader, we have stated the main results below as Theorem~\ref{thmE8}.\ 
The result is probably well-known to experts.

\bigskip

\noindent{\bf Supplemental material}

Supplemental data associated with this article can be found, in the electronic version, at 
%arXiv:1505.06420.
\hbox{http://arxiv.org/abs/1505.06420.}

%%%%%%%%%%%%%%%%%%%%%%%%%%%%%%%%%%%%%%%%%%%%%%%%%%%%%%%%%%%%%%%%%%%%%%%%%%%%%%%%%%%%%%%%%%%%%%%%%%%%%%%%%%%%%%%%%

\section{Integral lattices and their automorphism groups}\label{basicstuff}

We introduce some notation related to integral lattices 
and their automorphism groups and record some results that we will need. 

\smallskip

A \emph{lattice\/} $L$ is a finitely generated free $\Z$-module together with a rational-valued symmetric bilinear form $(\,.\,,\,.\,)$.\
All lattices in this note are assumed to be positive-definite.\
We let $O(L):={\rm Aut}(L)$ be the group of automorphisms (or \emph{isometries\/}) of $L$ \emph{considered as lattice\/}, i.e., 
the set  of automorphisms of the group $L$ that preserve the bilinear form.\ 
It is finite because of the assumed positive-definiteness of the bilinear form.\
The lattice $L$ is \emph{integral\/} if the bilinear form takes values in $\Z$,
and \emph{even\/} if the \emph{norm\/} $(x,\,x)$ belongs to $2\Z$ for all $x\in L$.\ An even lattice is
necessarily integral.\  

A \emph{finite quadratic space\/} $A=(A, q)$ is a finite abelian group $A$ equipped with a 
quadratic form $q: A\longrightarrow \Q/2\Z$.\ We denote the corresponding orthogonal group by $O(A)$.\ This is the subgroup
of ${\rm Aut}(A)$ that leaves $q$ invariant.

The \emph{dual lattice\/} of an integral lattice $L$ is 
$$L^*:=\{x\in L\otimes\Q \mid (x,y)\in \Z \mbox{\ for all\ } y\in L\}.$$ 
The {\it discriminant group\/} $L^*/L$ of an even lattice $L$ is equipped with the
{\it discriminant form\/} $q_L: L^*/L \rightarrow \Q/2\Z$, $x+L \mapsto (x,\,x)\ ({\rm mod}\ 2\Z)$.\ 
This turns $L^*/L$ into a finite quadratic space, called the {\it discriminant space\/} of $L$ and denoted 
$A_L:=(L^*/L, q_L)$. \ There is a natural induced action of $O(L)$ on $A_L$,  leading to a short exact sequence
$$1\longrightarrow O_0(L)   \longrightarrow O(L) \longrightarrow \overline{O}(L)  \longrightarrow 1,$$
where $\overline{O}(L)$ is the subgroup of $O(A_L)$ induced by $O(L)$ and $O_0(L)$ consists of the automorphisms
of $L$ which act trivially on $A_L$.

\smallskip

A sublattice  $K\subseteq L$ is called \emph{primitive\/} (in $L$) if $L/K$ is a free abelian group.\ 
We set 
$$K^\perp:=\{x\in L\mid (x,y)=0 \hbox{\ for all\ } y\in K\}.$$

Assume now that $L$ is even and \emph{unimodular,\/} i.e., $L^*=L$.\ 
If $K$ is primitive then there is an isomorphism of groups $i: A_K\stackrel{\cong}{\longrightarrow} A_{K^{\perp}}$ 
such that $q_{K^\perp}(i(a))=-q_{K}(a)$ for $a\in  A_K$.\ We can recover
$L$  from $K\oplus K^\perp$ by adjoining the cosets 
$$C:=\{(a,i(a))\mid a \in A_K\}\subseteq A_{K} \oplus   A_{K^\perp}.$$
See~\cite{Nikulin} for further details.\
The following is a special case of another result
(Propositions~1.4.1 and~1.6.1, loc.\ cit). 
\begin{prop}\label{gluing}
The equivalence classes of extensions of $K\oplus K^\perp$ to an even unimodular lattice $N$
with $K$ primitively embedded into $N$ are in bijective correspondence with double cosets
$ \overline{O}(K)\backslash O(A_K)/i^*\overline{O}(K^\perp)$, where 
$i^*: \overline{O}(K^\perp) \longrightarrow {O}(A_K)$ is defined by $g\mapsto  i^{-1} \circ g\circ i$. \qed
\end{prop}

\smallskip

Suppose that  $G\subseteq O(L)$ is a group of automorphisms of a lattice $L$.\
The \emph{invariant\/} and \emph{coinvariant\/} lattices for $G$ are 
\begin{eqnarray*}
L^G &=& \{x\in L\mid gx=x \hbox{ for all\ } g\in G\}, \\
L_G &=& (L^G)^\perp\  =\  \{x\in L\mid (x,y)=0 \hbox{ for all\ } y\in L^G\}
\end{eqnarray*}
respectively.\ They are both primitive sublattices of $L$.\
The restriction of the $G$-action to $L_G$ induces an {\it embedding\/}
$G\subseteq O(L_G)$.

If $G\subseteq O(L)$, we denote by $\widetilde{G}$ the pointwise
stabilizer of $L^G$ in $O(L)$.\ We always have $G\subseteq \widetilde{G}$
and $L^G = L^{\widetilde{G}}$.\ Moreover $N_{O(L)}(\widetilde{G})$ is 
the setwise stabilizer of $L^G$, and $N_{O(L)}(\widetilde{G})/\widetilde{G}$ is a (faithful) group of isometries of $L^G$.
\begin{lem}\label{lemO0}
Suppose that $L$ is even and unimodular.\ Then $\widetilde{G} \cong O_0(L_G)$.
\end{lem}
\pf As explained above, 
$L$ is obtained from  $L^G\oplus L_G$ by adjoining cosets 
$C:=\{(a,i(a)) \mid a \in A_{L^G}\}\subseteq A_{L^G} \oplus A_{L_G}$.\ 
Furthermore, in this case
$\widetilde{G}$ necessarily acts trivially on $A_{L_G}$, so that
$\widetilde{G}\subseteq O_0(L_G)$.

On the other hand, we can extend the $O_0(L_G)$-action on $L_G$ to a
trivial action on $L^G$.\ Since $O_0(L_G)$ acts trivially on $A_{L^G}\oplus A_{L_G}$,
the action on $L_G\oplus L^G$ extends to an action on $L$.\ 
The Lemma follows. \qquad \qed

\smallskip

A \emph{root\/} of $L$ is a primitive vector in $v\in L$ such that reflection in $(\Z v)^{\perp}$ is an isometry of $L$.
The \emph{root sublattice\/} of $L$ is the sublattice spanned by all roots.

We also note that the \emph{genus\/} of a positive-definite even lattice $L$ is determined by
the quadratic space $A_L$ together with the rank of $L$~\cite{Nikulin}.

\medskip

We recall the following fact:
\begin{lem}
A finite group $G$ has a unique minimal normal subgroup $N$ such that $G/N$
is a $2$-group.\ It is the subgroup generated by all elements of odd order. \hfill \qed
\end{lem}
We follow usual practice and set $N=O^2(G)$.

\begin{lem}\label{lemmaO2} Let $L$ be a lattice and assume that
$G\subseteq O(L)$ satisfies $G=\widetilde{G}$. Then
$$  O^2(G)\unlhd \widetilde{O^2(G)} \unlhd G.$$
\end{lem}
\pf Since $L^G \subseteq L^{O^2(G)}$ then  $\widetilde{O^2(G)} \subseteq \widetilde{G}=G$.\ 
Moreover, since $O^2(G) \unlhd G$ then $G$ acts on
$ L^{O^2(G)}$, and hence normalizes the pointwise stabilizer $\widetilde{O^2(G)}$ of this lattice.
\phantom{xxxxxxxxx}\qquad \qed

%%%%%%%%%%%%%%%%%%%%%%%%%%%%%%%%%%%%%%%%%%%%%%%%%%%%%%%%%%%%%%%%%%%%%%%%%%%%%%%

\section{Construction of the fixed-point lattices}\label{construction}

Recall~\cite{Conway2} that the $2^{24}$ cosets comprising $\Lambda/2\Lambda$ have representatives
$v$ which may be chosen to be \emph{short vectors\/}, i.e., 
$(v, v)\leq 8$.\ More precisely, if $(v, v)\leq 6$ then $\{v,-v\}$ are the \emph{only\/}
short representatives of $v+2\Lambda$; if $(v, v)=8$ then the short
vectors in $v+2\Lambda$ comprise a \emph{coordinate frame\/} $\{\pm w_1,\, \ldots,\, \pm w_{24}\}$,
where the $w_j$  are pairwise orthogonal vectors of norm~$8$.\ In particular, if $u\in\Lambda$ then
$u=v+2w$ for some $v$, $w\in\Lambda$ and $v$ a short vector, and if $(v, v)\leq 6$ then $v$ is 
unique up to sign.

It is well-known~\cite{Conway2} that $\Co$ acts \emph{transitively\/} on coordinate frames, the (setwise)
stabilizer of one such being the \emph{monomial group\/} $\mon$.

A sublattice $S\subseteq \Lambda$ is an \emph{$\S$-lattice\/} if, for every $u\in S$, the corresponding
short vector $v$ satisfies $(v, v)\leq 6$ and furthermore $w\in S$.\ This concept was introduced by 
Curtis~\cite{Curtis} who showed that there are exactly twelve isometry classes of $\S$-lattices.
The next result is a useful variant of a construction given in the Atlas~\cite{Atlas}.

\begin{prop}\label{lemmaodd} If $G=O^2(G)\subseteq \Co$, then one of the following holds:
\vspace{-2mm}
\begin{itemize}
\itemsep0em
\item[(a)] $\widetilde{G}\subseteq \mon$
\item[(b)] $\Lambda^G$ is an $\S$-lattice.
\end{itemize}
\end{prop}
\pf Let $u\in \Lambda^G$ with $u=v+2w$, where $v$, $w\in\Lambda$, and $(v, v)\leq 8$.\ 
Then $v+2\Lambda = u+2\Lambda$ is $G$-invariant.

First suppose that for every choice of $u$, we have $(v, v)\leq 6$.\ 
Then $\{\pm v\}$ are the only short vectors in $u+2\Lambda$, so this set is invariant under the
action of $G$.\ Then every odd order element in $G$ fixes $v$, and since $G=O^2(G)$ 
then $v\in\Lambda^G$.\ 
Then also $2w=u-v\in\Lambda^G$, and because $\Lambda^G$ is primitive then $w\in\Lambda^G$.\
So (b) holds in this case.

Otherwise, for some $u\in\Lambda^G$ we have $u=v+2w$ and $(v, v)=8$.\ 
Then because $\widetilde{G}$ fixes $u$, it acts on $u+2\Lambda$ and therefore stabilizes the unique coordinate frame
contained in this coset.\ So in this case (a) holds.  \qed

\begin{rem}\label{remS} For the stabilizer $G=\widetilde{G}$ of an $\S$-lattice one has always $G=O^2(G)$ and 
$G\not\subseteq\mon$.
\end{rem}
\pf This can easily be seen directly from the classification of $\S$-lattices and their stabilizers~\cite{Curtis} 
(cf.\ Table~\ref{maintable}).\ Note that $|G|$ does not divide $|\mon|$ so that part~(a) of 
Proposition~\ref{lemmaodd} fails by Lagrange's Theorem. \qed

\smallskip

Now assume that $G=\widetilde{G}\subseteq \Co$ with $L:=\Lambda^G\subseteq M:=\Lambda^{O^2(G)}$.\ 
By Lemma~\ref{lemmaO2} we have
$O^2(G)\unlhd \widetilde{O^2(G)}\unlhd G$, and $\widetilde{O^2(G)}$ is the pointwise
stabilizer of $M$.\ Thus  $L=M^G = M^{G/\widetilde{O^2(G)}}$ is the fixed-point sublattice of
$G/\widetilde{O^2(G)}$, which is a faithful $2$-group of isometries of $M$.\ 
Furthermore, by Proposition~\ref{lemmaodd} and Remark~\ref{remS}, either $\widetilde{O^2(G)}\subseteq \mon$ 
or $M$ is an $\S$-lattice.

\medskip 

This leads to the following general approach for finding all fixed-point lattices~$L$:
\vspace{-4mm}
\begin{itemize}
\itemsep0em
\item[(a)] Find all subgroups $H=O^2(H)\subseteq \mon$ and 
all pointwise stabilizers $H$ of $\S$-lattices (cf.\ Remark~\ref{remS}).
\item[(b)] For each such $H$, calculate $L=\Lambda^G$ where $H\unlhd G$ and $G/H$ is a $2$-group.
\end{itemize}

\medskip

We say that two pairs of lattices $(L_1,L_2)$ and $(L_1',L_2')$ are 
isometric if there are isometries $L_i\stackrel{\cong}{\longrightarrow}L_i'$
($i=1$, $2$).\

\smallskip

In order to make the enumeration of the fixed-point lattices outlined above effective, we 
iteratively compile a list of triples $(G,\Lambda^G,\Lambda_G)$ using the following procedure.

\begin{itemize}
\item[\bf Step 1:]
Select a representative $G$ from each conjugacy class of subgroups of $\mon$ 
satisfying $G=O^2(G)$.\ Construct $(G, \Lambda^G,\Lambda_G)$.
Select one triple for each isometry class of pairs $(\Lambda^G,\Lambda_G)$  of lattices, resulting
in a list of such triples.

\item[\bf Step 2:]
For each triple $(G,\Lambda^G,\Lambda_G)$, construct the pointwise-stabilizer 
$\widetilde{G}=O_0(\Lambda_G)$ of $\Lambda^G$ (cf.\ Lemma~\ref{lemO0}) in $\Co$ and replace $G$ by $\widetilde{G}$.

\item[\bf Step 3:]
For each triple $(G,\Lambda^G,\Lambda_G)$, compute the normalizer $N$ of $G$ in $\Co$. For each conjugacy
class $[g]$ in $N/G$, construct the group $H=\langle G, g \rangle$ and add
the triple $(H,\Lambda^{H},\Lambda_{H})$ to the list if $(\Lambda^{H},\Lambda_{H})$ is not isometric 
to a pair of lattices already present.

\item[\bf Step 4:] Repeat Steps 2 and 3 until the list is saturated.

\end{itemize}

This results in the list of $290$ triples which, along with accompanying data, are described in Table~\ref{maintable}.

\medskip

We explain now why the triples resulting from Steps 1--4 produce the desired list of 
orbits of fixed-point lattices, thereby proving Theorem~\ref{mainthm}.

\smallskip

First, notice that if $(G, \Lambda^G,\Lambda_G)$ and $(H, \Lambda^H, \Lambda_H)$ are 
distinct triples on the final list, then $(\Lambda^G,\Lambda_G)$ and $(\Lambda^H, \Lambda_H)$
are \emph{not\/} isometric.\ Therefore, $\Lambda^G$ and $\Lambda^H$ certainly lie in distinct
$\Co$-orbits, since an element of $\Co$ mapping $\Lambda^G$ onto $\Lambda^H$ is an isometry
that also induces an isometry of $\Lambda_G$ onto $\Lambda_H$.

\smallskip

Next we show that every $\Co$-orbit of fixed-point lattices has a representative that occurs 
in a triple on the final list.\ First we verify that the isometry classes
$(\Lambda^G,\Lambda_G)$ already determine the orbits of fixed-point lattices. 
\begin{prop}\label{uniqueness}
For each entry in Table~\ref{maintable}, the isometry class of the pair $(\Lambda^G,\Lambda_G)$
uniquely determines the $\Co$-orbit of $\Lambda^G$.
\end{prop}

\pf For each pair $(\Lambda^G,\Lambda_G)$,
we determine all isomorphism classes of extensions of $\Lambda^G\oplus \Lambda_G$ to an even unimodular
lattice $N$ (i.e.~the even unimodular overlattices of $\Lambda^G\oplus \Lambda_G$) 
by computing the double cosets for $\overline{O}(\Lambda^G)\times i^*(\overline{O}(\Lambda^G))$ in $O(A_{\Lambda^G})$ 
(cf.\ Proposition~\ref{gluing}).\
Among the resulting lattices $N$, it turns out there is always \emph{exactly one\/}
equivalence class with minimal norm $4$, so that it must be isometric to $\Lambda$.

It follows that $(\Lambda^G,\Lambda_G)$ uniquely determines the $\Co$-orbit of $\Lambda^G$ since
two extensions $L$ and $L'$ of $\Lambda^G\oplus \Lambda_G$ are by definition 
equivalent if there is an isometry between $L$ and $L'$ which stabilizes  $\Lambda^G\oplus \Lambda_G$ setwise,
i.e.\ after identifications of $L$ and $L'$ with $\Lambda$, 
the corresponding sublattices $\Lambda^G\oplus \Lambda_G$ of $\Lambda$ can be mapped to
each other by an element of $\Co$.  \qed

\smallskip

Next, all $\S$-lattices and their stabilizers appear in Table~\ref{maintable}.\ Indeed,  
the twelve lattices  $\Lambda^G$ numbered
$35$, $101$, $122$, $163$, $167$, $222$, $223$, $225$, $230$, $273$, $274$ and
$290$ have the two properties: $G=O^2(G)$ and $|G|$ does \emph{not\/} divide $|\mon|$.\ 
By Proposition~\ref{lemmaodd},  each $\Lambda^G$ is an $\S$-lattice.\ 
According to Curtis~\cite{Curtis} there are exactly twelve $\Co$-orbits of $\S$-lattices, 
so indeed they all appear in Table~\ref{maintable}.

Along with the $\S$-lattices, Step~3 ensures that with a fixed-point lattice $\Lambda^H$, 
all fixed-point lattices $\Lambda^G$ also occur in a triple whenever $H\unlhd G$ and
$G/H$ is a $2$-group as the following proposition shows.
\begin{prop}\label{2groups}
Assume that for $G\subseteq \Co$ the triple  $(\widetilde{O^2(G)},\Lambda^{\widetilde{O^2(G)}},\Lambda_{\widetilde{O^2(G)}})$
is contained on the list in Table~\ref{maintable}.\ Then $(\widetilde{G}, \Lambda^{\widetilde{G}},\Lambda_{\widetilde{G}})$ is also contained in the list.
\end{prop}
\pf Because it is a $2$-group, $G/O^2(G)$ has  a central series
$$O^2(G)=H_0 \unlhd H_1 \unlhd \cdots \unlhd H_n =G$$
with each $H_i\unlhd G$ and $|H_{i+1}/H_i|=2$, and
$$ \Lambda^{O^2(G)}=\Lambda^{H_0} \supseteq \Lambda^{H_1}  \supseteq  \cdots   \supseteq \Lambda^{H_n} =\Lambda^G.$$
$G$ acts on each $\Lambda^{H_i}$, and hence  normalizes $\widetilde{H_i}$.\ Using $H_{i+1}/(H_{i+1}\cap \widetilde{H_i})\cong H_{i+1}\widetilde{H_i}/\widetilde{H_i}$ and
$H_{i}\subseteq H_{i+1}\cap \widetilde{H_i}$ we conclude that  $[H_{i+1}\widetilde{H_i}:\widetilde{H_i}]  \leq 2$.\
Thus Steps~2 and~3 guarantee that $(\widetilde{H_{i+1}},\Lambda^{\widetilde{H_{i+1}}}, \Lambda_{\widetilde{H_{i+1}}})$ 
is on the list whenever $(\widetilde{H_i},\Lambda^{\widetilde{H_i}},\Lambda_{\widetilde{H_i}})$ is.
 
Since, by assumption, the triple $(\widetilde{H_0},\Lambda^{\widetilde{H_0}},\Lambda_{\widetilde{H_0}})$ 
is contained in the list, it follows inductively that  
$(\widetilde{G},\Lambda^{\widetilde{G}},\Lambda_{\widetilde{G}})$ is too. \qed

\smallskip

Together with the results of the computation, we have established Theorem~\ref{mainthm}.

\medskip

We describe now some more details for the implementation of Steps 1 to 3 with the computer algebra
system MAGMA.

\smallskip

We realized the Conway group $\Co$ as a matrix group of integral $24\times 24$-matrices and as a permutation
group on the $196,560$ vectors of norm~$4$. We also determined an explicit isomorphism which allows us to 
evaluate a computation in the most appropriate realization. 

For Step~1, we started with the list of conjugacy classes of non-$2$-groups inside $2^{12}{:}M_{24}$.
In~\cite{K32} we had already shown:
\begin{thm}\label{non2group} With respect to conjugation in $\mon$, 
there are $279,343$ conjugacy classes of subgroups of $\mon$ which are \emph{not\/} $2$-groups.
\end{thm}
From these classes  we selected those groups $G$ which satisfy $G=O^2(G)$. This was done by computing 
$O^2(G)$ as the normal subgroup of $G$ generated by $p$-Sylow subgroups for all $p\not=2$.  
This resulted in a list of $3755$ groups.  
For these groups we computed the pairs $(\Lambda^G,\Lambda_G)$ of sublattices inside $\Lambda$
and checked for isometry by the implemented lattice functions in MAGMA.

For Step~2, we can compute $\widetilde{G}$ abstractly as the group $O_{0}(\Lambda_G)$. However, to 
realize $\widetilde{G}$ as a subgroup of $\Co$ we realized in addition $\Co$ 
as a matrix group over the finite field $\F_2$ acting on $\Lambda/2\Lambda$. 
This allowed us to compute that stabilizer of $\Lambda^G/2\Lambda$ in $\Co$.

Step 3 can easily be done by the implemented group theory functions in MAGMA.

\bigskip

\paragraph{Remarks on the $E_8$-root lattice.} The $E_8$-root lattice is the unique even, unimodular,
positive-definite lattice of rank~$8$ and its automorphism group is the corresponding Weyl group.\
The problem of determining the orbits of fixed-point sublattices and stabilizer groups for this 
lattice and its automorphism group also has some interest attached to it.\

It follows from a Theorem of Steinberg~(\cite{Steinberg}, Thm.~15) that the stabilizer 
of a sublattice of a root lattice inside the corresponding Weyl group is a reflection group. 
The conjugacy classes of reflections subgroups for $W(E_8)$ are known, cf.~\cite{DPR}, Table~5.

We desist from further discussion, contenting ourselves with the statement of the result, which must be well-known 
to experts.
\begin{thm}\label{thmE8}
In its action on the $E_8$-root lattice, the Weyl group of type $E_8$ has $41$ orbits of fixed-point sublattices.\ 
These are in bijective correspondence with the isomorphism types
of full subgraphs of the Coxeter graph for $E_8$, the lattice-stabilizers being the Coxeter groups
determined by these subgraphs.\ The coinvariant lattices are the corresponding root lattices.
\qed
\end{thm}

%%%%%%%%%%%%%%%%%%%%%%%%%%%%%%%%%%%%%%%%%%%%%%%%%%%%%%%%%%%%%%%%%%%%%%%%%%%%%%%%%%%%%%%%%%%%%%%%%%%%%%%%%%%%%%%%%
 
\section{The $290$ fixed-point lattices}\label{SSA}

This section describes the fixed-point lattices and discusses several observations regarding the resulting tables. 
Table~\ref{maintable} provides information about the $290$ orbits
of fixed-point lattices $L=\Lambda^G$ inside $\Lambda$.\ 
For a given $L$, the group $G$ listed is the full pointwise
stabilizer $\Co$, i.e., $G=\widetilde{G}=O_0(\Lambda_G)$.\

In addition, we provide the following {\it electronic supplementary material.\/} 
Supplementary Table~2 consists of the Gram matrices of each $\Lambda^G$.
Supplementary Table~3 gives partial information about the lattice structure of the $290$ orbits. 
We have also added a text file in Magma format which containes a coinvariant lattice from each orbit, the 
corresponding fixed-point lattice and generator matrices for its stabilizers.

\medskip

%%%%%%%%%%%%%%%%%%%%%%%%%%%%%%%%%%%%%%%%%%%%%%%%%%%%%%%%%

{\LTcapwidth=\linewidth
\small\noindent
\begin{longtable}{rrrllrrrrrrl}
\caption{Orbits of fixed-point lattices.
The columns provide the following information:
number of $\Lambda^G$ (no.);
rank of $\Lambda^G$ (rk);
order of $G$ (order).\
Information about the group structure of $G$ ($G$).\ Here, $[n]$ denotes an unspecified group of order $n$ and $p^n$ an elementary abelian group of the same order.\
Sometimes we list the standard name for the group or the number of $G$ in the database of small groups.\
The genus symbol for $\Lambda^G$ without the signature information (genus);
rank of $\Lambda^G$ minus the rank of $A_{\Lambda^G}$ ($\alpha$);
%\item[$O(A)$:] Order of orthogonal group of discriminant group $A_{\Lambda^G}$;
index of $\overline{O}(\Lambda_G)$ in $O(A_{\Lambda^G})$ ($\overline{i}_G$);
index of $\overline{O}(\Lambda^G)$ in $O(A_{\Lambda^G})$ ($\overline{i}^G$);
index of ${\rm N}_{{\rm Co}_0}(G)/G$ in $O(L^G)$ (ind); 
%\item[$h_G$:] number of lattices in the genus of $\Lambda_G$.
number of lattices in the genus of $\Lambda^G$  ($h^G$);
number of Niemeier lattices with roots into which $\Lambda_G$ embeds ($N$);
case type ([$M_{23}$]: $G\subseteq M_{23}$, [$M_{24}$]: $G\subseteq M_{24}$ and not [$M_{23}$];
[${\rm Mon}_a$]: $G\subseteq \mon$ but not [$M_{23}$], [$M_{24}$] and $G=T{:}H$ with 
$H\subseteq M_{24}$ and $T=G\cap 2^{12}$;
[${\rm Mon}_b$]: $G\subseteq \mon$ but not [$M_{23}$], [$M_{24}$], [${\rm Mon}_a$];
[$-$]: $G\subsetneq\mon$ but not [S], [S]: $|G|\nmid |\mon|$;
[*]: $\widetilde{O^2(G)}=G$) (type).
%An empty entry means that quantity is not determined.
}
\label{maintable}  \\ 
\mbox{no.} & rk & $\!\!$\mbox{order} & $G$ & $\!\!$genus & $\!\!\alpha$ &  $\overline{i}_G$ & $\overline{i}^G$ & ind  & $h^G$ & $N$  & type \\ \hline
\endfirsthead 
\caption[]{Orbits of fixed-point lattices}\\
\mbox{no.} & rk & $\!\!$\mbox{order} & $G$ & $\!\!$genus & $\!\!\alpha$ &  $\overline{i}_G$ & $\overline{i}^G$ & ind  & $h^G$ & $N$ & type  \\ \hline
\endhead

1 & 
24 & 
1 & 
$1$& 
$1$ & 
24 & 
1 & 
1 & 
1 & 
24 & 
23 & 
${M}_{23}$* 
\\  
 
\hline  
%%%%%%%%%%%%%%%%%%%%%%%%%%%%%%%%%%%%%%%%%%%%%%%%%%%%%%%%%%%%%%%%%%%%%% 

2 & 
16 & 
2 & 
$2$& 
$2_{\rm I\!I}^{+8}$ & 
8 & 
1 & 
2 & 
1 & 
24 & 
17 & 
${M}_{23}$ 
\\  
 
\hline  
%%%%%%%%%%%%%%%%%%%%%%%%%%%%%%%%%%%%%%%%%%%%%%%%%%%%%%%%%%%%%%%%%%%%%% 

3 & 
12 & 
4 & 
$2^2$& 
$2_{\rm I\!I}^{-6}4_{\rm I\!I}^{-2}$ & 
4 & 
40 & 
40 & 
12 & 
7 & 
7 & 
$M_{23}$ 
\\  
 
4 & 
12 & 
3 & 
$3$& 
$3^{+6}$ & 
6 & 
1 & 
1 & 
1 & 
10 & 
8 & 
$M_{23}$* 
\\  
 
5 & 
12 & 
2 & 
$2$& 
$2_{4}^{+12}$ & 
0 & 
104448 & 
104448 & 
5040 & 
3 & 
11 & 
${M}_{24}$ 
\\  
 
\hline  
%%%%%%%%%%%%%%%%%%%%%%%%%%%%%%%%%%%%%%%%%%%%%%%%%%%%%%%%%%%%%%%%%%%%%% 
6 & 
10 & 
8 & 
$2^3$& 
$2_{\rm I\!I}^{+6}4_{2}^{+2}$ & 
2 & 
135 & 
36 & 
30 & 
4 & 
2 & 
$M_{23}$ 
\\  
 
7 & 
10 & 
6 & 
$S_3$& 
$2_{\rm I\!I}^{-2}3^{+5}$ & 
5 & 
1 & 
2 & 
1 & 
13 & 
7 & 
$M_{23}$ 
\\  
 
8 & 
10 & 
4 & 
$2^2$& 
$2_{\rm I\!I}^{+8}4_{2}^{+2}$ & 
0 & 
45696 & 
26112 & 
2520 & 
2 & 
6 & 
$M_{24}$ 
\\  
 
9 & 
10 & 
4 & 
$4$& 
$2_{2}^{+2}4_{\rm I\!I}^{+4}$ & 
4 & 
1 & 
2 & 
1 & 
8 & 
7 & 
$M_{23}$ 
\\  
 
\hline  
%%%%%%%%%%%%%%%%%%%%%%%%%%%%%%%%%%%%%%%%%%%%%%%%%%%%%%%%%%%%%%%%%%%%%% 
10 & 
9 & 
16 & 
$2^4$& 
$2_{\rm I\!I}^{+6}8_{1}^{+1}$ & 
2 & 
2 & 
1 & 
2 & 
3 & 
1 & 
$M_{23}$ 
\\  
 
11 & 
9 & 
16 & 
$2^4$& 
$2_{\rm I\!I}^{+8}4_{1}^{+1}$ & 
0 & 
2295 & 
136 & 
270 & 
1 & 
1 & 
$M_{24}$ 
\\  
 
12 & 
9 & 
8 & 
$2^3$& 
$2_{\rm I\!I}^{+8}8_{1}^{+1}$ & 
0 & 
11200 & 
960 & 
840 & 
2 & 
4 & 
$M_{24}$ 
\\  
 
13 & 
9 & 
8 & 
$D_8$& 
$4_{1}^{+5}$ & 
4 & 
1 & 
4 & 
1 & 
8 & 
5 & 
$M_{23}$ 
\\  
 
\hline  
%%%%%%%%%%%%%%%%%%%%%%%%%%%%%%%%%%%%%%%%%%%%%%%%%%%%%%%%%%%%%%%%%%%%%% 
14 & 
8 & 
512 & 
$[2^9]$& 
$2_{\rm I\!I}^{+8}$ & 
0 & 
2 & 
1 & 
2 & 
1 & 
- & 
${\rm Mon}_a$ 
\\  
 
15 & 
8 & 
18 & 
$A_{3,3}$& 
$3^{+4}9^{-1}$ & 
3 & 
3 & 
3 & 
2 & 
3 & 
3 & 
$M_{23}$* 
\\  
 
16 & 
8 & 
16 & 
$2^4$& 
$2_{\rm I\!I}^{+6}4_{\rm I\!I}^{+2}$ & 
0 & 
840 & 
64 & 
105 & 
1 & 
2 & 
$M_{24}$ 
\\  
 
17 & 
8 & 
16 & 
$\Gamma_2a_1$& 
$2_{\rm I\!I}^{+2}4_{0}^{+4}$ & 
2 & 
6 & 
16 & 
3 & 
3 & 
2 & 
$M_{23}$ 
\\  
 
18 & 
8 & 
12 & 
$D_{12}$& 
$2_{\rm I\!I}^{+4}3^{+4}$ & 
4 & 
1 & 
12 & 
1 & 
8 & 
6 & 
$M_{23}$ 
\\  
 
19 & 
8 & 
12 & 
$A_4$& 
$2_{\rm I\!I}^{-2}4_{\rm I\!I}^{-2}3^{+2}$ & 
4 & 
1 & 
1 & 
1 & 
9 & 
4 & 
$M_{23}$* 
\\  
 
20 & 
8 & 
10 & 
$D_{10}$& 
$5^{+4}$ & 
4 & 
1 & 
1 & 
1 & 
5 & 
4 & 
$M_{23}$* 
\\  
 
21 & 
8 & 
8 & 
$[2^3]$\ (\#3)& 
$2_{\rm I\!I}^{+4}4_{\rm I\!I}^{+4}$ & 
0 & 
512 & 
8192 & 
36 & 
1 & 
6 & 
$M_{24}$ 
\\  
 
22 & 
8 & 
6 & 
$S_3$& 
$3^{+8}$ & 
0 & 
56862 & 
56862 & 
1920 & 
1 & 
7 & 
$M_{24}$* 
\\  
 
23 & 
8 & 
4 & 
$2^2$& 
$2_{4}^{+4}4_{4}^{+4}$ & 
0 & 
16 & 
1024 & 
6 & 
3 & 
5 & 
$M_{24}$ 
\\  
 
\hline  
%%%%%%%%%%%%%%%%%%%%%%%%%%%%%%%%%%%%%%%%%%%%%%%%%%%%%%%%%%%%%%%%%%%%%% 
24 & 
7 & 
48 & 
$2^4{:}3$& 
$2_{\rm I\!I}^{-4}8_{1}^{+1}3^{-1}$ & 
2 & 
1 & 
1 & 
1 & 
4 & 
1 & 
$M_{23}$* 
\\  
 
25 & 
7 & 
32 & 
$2^5$& 
$2_{\rm I\!I}^{+6}8_{7}^{+1}$ & 
0 & 
28 & 
8 & 
7 & 
1 & 
1 & 
$M_{24}$ 
\\  
 
26 & 
7 & 
32 & 
$\Gamma_4a_1$& 
$2_{\rm I\!I}^{+2}4_{6}^{+2}8_{1}^{+1}$ & 
2 & 
1 & 
2 & 
1 & 
4 & 
1 & 
$M_{23}$ 
\\  
 
27 & 
7 & 
32 & 
$[2^5]$\ (\#46)& 
$2_{\rm I\!I}^{+4}4_{7}^{+3}$ & 
0 & 
60 & 
64 & 
15 & 
1 & 
1 & 
$M_{24}$ 
\\  
 
28 & 
7 & 
32 & 
$\Gamma_5a_1$& 
$4_{7}^{+5}$ & 
2 & 
1 & 
16 & 
1 & 
2 & 
1 & 
$M_{23}$ 
\\  
 
29 & 
7 & 
24 & 
$S_4$& 
$4_{3}^{+3}3^{+2}$ & 
4 & 
1 & 
2 & 
1 & 
11 & 
4 & 
$M_{23}$ 
\\  
 
30 & 
7 & 
8 & 
$[2^3]$\ (\#3)& 
$2_{6}^{+4}4_{\rm I\!I}^{+2}8_{1}^{+1}$ & 
0 & 
16 & 
128 & 
6 & 
2 & 
3 & 
$M_{24}$ 
\\  
 
31 & 
7 & 
8 & 
$Q_8$& 
$2_{3}^{+3}8_{\rm I\!I}^{-2}$ & 
2 & 
1 & 
1 & 
1 & 
4 & 
3 & 
$M_{23}$ 
\\  
 
32 & 
7 & 
8 & 
$2^3$& 
$2_{\rm I\!I}^{+2}4_{7}^{+5}$ & 
0 & 
32 & 
4096 & 
6 & 
1 & 
3 & 
$M_{24}$ 
\\  
 
\hline  
%%%%%%%%%%%%%%%%%%%%%%%%%%%%%%%%%%%%%%%%%%%%%%%%%%%%%%%%%%%%%%%%%%%%%% 
33 & 
6 & 
1536 & 
$[2^93]$& 
$2_{\rm I\!I}^{-6}3^{-1}$ & 
0 & 
1 & 
1 & 
1 & 
1 & 
- & 
${\rm Mon}_a$* 
\\  
 
34 & 
6 & 
1024 & 
$[2^{10}]$& 
$2_{\rm I\!I}^{+4}4_{6}^{+2}$ & 
0 & 
1 & 
2 & 
1 & 
1 & 
- & 
${\rm Mon}_a$ 
\\  
 
35 & 
6 & 
486 & 
$3^{1+4}{:2}$& 
$3^{+5}$ & 
1 & 
1 & 
1 & 
1 & 
1 & 
- & 
S* 
\\  
 
36 & 
6 & 
192 & 
$[2^63]$\ (\#1541)& 
$2_{\rm I\!I}^{+4}4_{6}^{+2}$ & 
0 & 
2 & 
2 & 
1 & 
1 & 
1 & 
$M_{24}$* 
\\  
 
37 & 
6 & 
192 & 
$4^2.A_4$& 
$2_{\rm I\!I}^{-2}8_{6}^{-2}$ & 
2 & 
1 & 
1 & 
1 & 
2 & 
1 & 
$M_{23}$* 
\\  
 
38 & 
6 & 
96 & 
$2^4.D_6$& 
$2_{\rm I\!I}^{-2}4_{7}^{+1}8_{1}^{+1}3^{-1}$ & 
2 & 
1 & 
2 & 
1 & 
5 & 
1 & 
$M_{23}$ 
\\  
 
39 & 
6 & 
72 & 
$A_{4,3}$& 
$4_{\rm I\!I}^{-2}3^{-3}$ & 
3 & 
1 & 
2 & 
1 & 
3 & 
2 & 
$M_{23}$* 
\\  
 
40 & 
6 & 
64 & 
$[2^6]$\ (\#264)& 
$2_{\rm I\!I}^{+2}4_{6}^{+4}$ & 
0 & 
6 & 
64 & 
3 & 
1 & 
1 & 
$M_{24}$ 
\\  
 
41 & 
6 & 
64 & 
$[2^6]$\ (\#266)& 
$2_{6}^{+2}4_{\rm I\!I}^{+4}$ & 
0 & 
1 & 
32 & 
1 & 
1 & 
- & 
${\rm Mon}_b$ 
\\  
 
42 & 
6 & 
64 & 
$[2^6]$\ (\#202)& 
$2_{\rm I\!I}^{+4}4_{5}^{+1}8_{1}^{+1}$ & 
0 & 
15 & 
12 & 
5 & 
1 & 
1 & 
$M_{24}$ 
\\  
 
43 & 
6 & 
64 & 
$\Gamma_{25}a_1$& 
$4_{5}^{+3}8_{1}^{+1}$ & 
2 & 
1 & 
4 & 
1 & 
3 & 
1 & 
$M_{23}$ 
\\  
 
44 & 
6 & 
60 & 
$A_5$& 
$2_{\rm I\!I}^{-2}3^{+1}5^{-2}$ & 
4 & 
1 & 
1 & 
1 & 
6 & 
3 & 
$M_{23}$* 
\\  
 
45 & 
6 & 
48 & 
$2\times S_4$& 
$2_{\rm I\!I}^{+2}4_{2}^{+2}3^{+2}$ & 
2 & 
3 & 
6 & 
2 & 
5 & 
2 & 
$M_{23}$ 
\\  
 
46 & 
6 & 
36 & 
$3^2.4$& 
$2_{2}^{+2}3^{+2}9^{-1}$ & 
3 & 
1 & 
1 & 
1 & 
4 & 
2 & 
$M_{23}$ 
\\  
 
47 & 
6 & 
36 & 
$S_{3,3}$& 
$2_{\rm I\!I}^{-2}3^{+3}9^{-1}$ & 
2 & 
3 & 
18 & 
2 & 
3 & 
3 & 
$M_{23}$ 
\\  
 
48 & 
6 & 
32 & 
$\Gamma_5a_2$& 
$2_{6}^{+2}4_{\rm I\!I}^{+4}$ & 
0 & 
30 & 
32 & 
15 & 
1 & 
1 & 
$M_{24}$ 
\\  
 
49 & 
6 & 
32 & 
$[2^5]$\ (\#27)& 
$2_{\rm I\!I}^{+2}4_{6}^{+4}$ & 
0 & 
36 & 
64 & 
9 & 
1 & 
2 & 
$M_{24}$ 
\\  
 
50 & 
6 & 
24 & 
$[2^33]$\ (\#14)& 
$2_{\rm I\!I}^{-6}3^{+3}$ & 
0 & 
40 & 
240 & 
12 & 
1 & 
2 & 
$M_{24}$ 
\\  
 
51 & 
6 & 
24 & 
$S_4$& 
$2_{\rm I\!I}^{+4}4_{\rm I\!I}^{-2}3^{+1}$ & 
0 & 
32 & 
32 & 
6 & 
1 & 
3 & 
$M_{24}$ 
\\  
 
52 & 
6 & 
21 & 
$F_{21}$& 
$7^{+3}$ & 
3 & 
1 & 
1 & 
1 & 
3 & 
2 & 
$M_{23}$* 
\\  
 
53 & 
6 & 
20 & 
${\rm Hol}(4)$& 
$2_{2}^{+2}5^{+3}$ & 
3 & 
1 & 
2 & 
1 & 
6 & 
4 & 
$M_{23}$ 
\\  
 
54 & 
6 & 
16 & 
$ [2^4]$\ (\#11)& 
$2_{\rm I\!I}^{+2}4_{5}^{+3}8_{1}^{+1}$ & 
0 & 
16 & 
256 & 
6 & 
2 & 
3 & 
$M_{24}$ 
\\  
 
55 & 
6 & 
16 & 
$\Gamma_3a_2$& 
$2_{5}^{+1}4_{1}^{+1}8_{\rm I\!I}^{+2}$ & 
2 & 
1 & 
2 & 
1 & 
4 & 
3 & 
$M_{23}$ 
\\  
 
56 & 
6 & 
12 & 
$A_4$& 
$2_{6}^{+2}4_{\rm I\!I}^{+4}$ & 
0 & 
32 & 
32 & 
6 & 
1 & 
3 & 
$M_{24}$* 
\\  
 
57 & 
6 & 
8 & 
$2^3$& 
$4_{6}^{+6}$ & 
0 & 
60 & 
2048 & 
15 & 
1 & 
3 & 
$M_{24}$ 
\\  
 
58 & 
6 & 
8 & 
$2^3$& 
$2_{2}^{+2}4_{3}^{+3}8_{1}^{+1}$ & 
0 & 
1 & 
64 & 
1 & 
2 & 
1 & 
${\rm Mon}_a$ 
\\  
 
59 & 
6 & 
8 & 
$2\times 4$& 
$2_{4}^{+4}8_{2}^{+2}$ & 
0 & 
6 & 
16 & 
3 & 
2 & 
3 & 
$M_{24}$ 
\\  
 
60 & 
6 & 
6 & 
$S_3$& 
$2_{6}^{+6}3^{-2}$ & 
0 & 
1 & 
2 & 
1 & 
2 & 
1 & 
${\rm Mon}_b$ 
\\  
 
61 & 
6 & 
6 & 
$S_3$& 
$2_{2}^{+6}3^{-4}$ & 
0 & 
1 & 
720 & 
1 & 
2 & 
1 & 
${\rm Mon}_b$ 
\\  
 
62 & 
6 & 
6 & 
$6$& 
$2_{\rm I\!I}^{-6}3^{-5}$ & 
0 & 
1 & 
51840 & 
1 & 
1 & 
- & 
${\rm Mon}_a$ 
\\  
 
63 & 
6 & 
6 & 
$6$& 
$2_{4}^{-6}3^{-3}$ & 
0 & 
2 & 
48 & 
1 & 
3 & 
3 & 
${\rm Mon}_a$ 
\\  
 
64 & 
6 & 
4 & 
$4$& 
$4_{6}^{+6}$ & 
0 & 
64 & 
2048 & 
6 & 
1 & 
7 & 
$M_{24}$ 
\\  
 
\hline  
%%%%%%%%%%%%%%%%%%%%%%%%%%%%%%%%%%%%%%%%%%%%%%%%%%%%%%%%%%%%%%%%%%%%%% 
65 & 
5 & 
6144 & 
$[2^{11}3]$& 
$2_{\rm I\!I}^{-4}8_{5}^{-1}$ & 
0 & 
1 & 
1 & 
1 & 
1 & 
- & 
${\rm Mon}_a$* 
\\  
 
66 & 
5 & 
3072 & 
$[2^93^2]$& 
$2_{\rm I\!I}^{-4}4_{7}^{+1}3^{-1}$ & 
0 & 
1 & 
2 & 
1 & 
1 & 
- & 
${\rm Mon}_a$ 
\\  
 
67 & 
5 & 
2048 & 
$ [2^{11}]$& 
$2_{\rm I\!I}^{+2}4_{5}^{+3}$ & 
0 & 
1 & 
4 & 
1 & 
1 & 
- & 
${\rm Mon}_a$ 
\\  
 
68 & 
5 & 
972 & 
$3^{1+4}{:}2.2$& 
$2_{1}^{+1}3^{-4}$ & 
1 & 
1 & 
1 & 
1 & 
1 & 
- & 
S 
\\  
 
69 & 
5 & 
960 & 
$M_{20}$& 
$2_{\rm I\!I}^{-2}8_{1}^{+1}5^{-1}$ & 
2 & 
1 & 
1 & 
1 & 
3 & 
1 & 
$M_{23}$* 
\\  
 
70 & 
5 & 
384 & 
$4^2.S_4$& 
$4_{3}^{+1}8_{2}^{+2}$ & 
2 & 
1 & 
2 & 
1 & 
2 & 
1 & 
$M_{23}$ 
\\  
 
71 & 
5 & 
384 & 
$[2^6].S_3$\ (\#20164)& 
$2_{\rm I\!I}^{+2}4_{5}^{+3}$ & 
0 & 
6 & 
4 & 
3 & 
1 & 
1 & 
$M_{24}$ 
\\  
 
72 & 
5 & 
360 & 
$A_6$& 
$4_{5}^{-1}3^{+2}5^{+1}$ & 
3 & 
1 & 
1 & 
1 & 
4 & 
2 & 
$M_{23}$* 
\\  
 
73 & 
5 & 
288 & 
$A_{4,4}$& 
$2_{\rm I\!I}^{+2}8_{1}^{+1}3^{+2}$ & 
2 & 
1 & 
2 & 
1 & 
3 & 
1 & 
$M_{23}$* 
\\  
 
74 & 
5 & 
192 & 
$2^4{:}D_{12}$& 
$4_{2}^{-2}8_{1}^{+1}3^{-1}$ & 
2 & 
1 & 
4 & 
1 & 
4 & 
1 & 
$M_{23}$ 
\\  
 
75 & 
5 & 
192 & 
$[2^63]$\ (\#1538)& 
$2_{\rm I\!I}^{-4}8_{7}^{+1}3^{-1}$ & 
0 & 
10 & 
10 & 
4 & 
1 & 
1 & 
$M_{24}$ 
\\  
 
76 & 
5 & 
192 & 
$(Q_8*Q_8){:}S_3$& 
$4_{7}^{-3}3^{+1}$ & 
2 & 
1 & 
1 & 
1 & 
2 & 
1 & 
$M_{23}$* 
\\  
 
77 & 
5 & 
168 & 
$ L_2(7)$& 
$4_{1}^{+1}7^{+2}$ & 
3 & 
1 & 
1 & 
1 & 
4 & 
2 & 
$M_{23}$* 
\\  
 
78 & 
5 & 
128 & 
$ [2^7]$\ (\#2326)& 
$4_{5}^{+5}$ & 
0 & 
1 & 
64 & 
1 & 
1 & 
- & 
${\rm Mon}_a$ 
\\  
 
79 & 
5 & 
128 & 
$[2^7]$\ (\#1758)& 
$2_{\rm I\!I}^{-2}4_{\rm I\!I}^{-2}8_{1}^{+1}$ & 
0 & 
6 & 
16 & 
3 & 
1 & 
1 & 
$M_{24}$ 
\\  
 
80 & 
5 & 
128 & 
$ [2^7]$\ (\#1759)& 
$2_{6}^{+2}4_{\rm I\!I}^{+2}8_{7}^{+1}$ & 
0 & 
1 & 
8 & 
1 & 
1 & 
- & 
${\rm Mon}_b$ 
\\  
 
81 & 
5 & 
128 & 
$ [2^7]$\ (\#1755)& 
$2_{\rm I\!I}^{+2}4_{4}^{+2}8_{1}^{+1}$ & 
0 & 
6 & 
16 & 
3 & 
1 & 
1 & 
$M_{24}$ 
\\  
 
82 & 
5 & 
120 & 
$S_5$& 
$4_{3}^{-1}3^{+1}5^{-2}$ & 
3 & 
1 & 
2 & 
1 & 
5 & 
3 & 
$M_{23}$ 
\\  
 
83 & 
5 & 
96 & 
$[2^5.3]$\ (\#226)& 
$2_{\rm I\!I}^{+4}4_{1}^{+1}3^{+2}$ & 
0 & 
15 & 
20 & 
6 & 
1 & 
1 & 
$M_{24}$ 
\\  
 
84 & 
5 & 
72 & 
$M_9$& 
$2_{3}^{+3}3^{-1}9^{-1}$ & 
2 & 
1 & 
1 & 
1 & 
4 & 
2 & 
$M_{23}$ 
\\  
 
85 & 
5 & 
72 & 
$3^2.D_8$& 
$4_{1}^{+1}3^{+2}9^{-1}$ & 
2 & 
1 & 
2 & 
1 & 
3 & 
2 & 
$M_{23}$ 
\\  
 
86 & 
5 & 
64 & 
$[2^6]$\ (\#226)& 
$4_{5}^{+5}$ & 
0 & 
30 & 
64 & 
15 & 
1 & 
1 & 
$M_{24}$ 
\\  
 
87 & 
5 & 
48 & 
$T_{48}$& 
$2_{7}^{+1}8_{\rm I\!I}^{-2}3^{-1}$ & 
2 & 
1 & 
2 & 
1 & 
4 & 
3 & 
$M_{23}$* 
\\  
 
88 & 
5 & 
48 & 
$[2^43]$\ (\#48)& 
$2_{\rm I\!I}^{-2}4_{3}^{+3}3^{+1}$ & 
0 & 
16 & 
32 & 
6 & 
2 & 
3 & 
$M_{24}$ 
\\  
 
89 & 
5 & 
32 & 
$[2^5]$\ (\#43)& 
$2_{5}^{+3}8_{\rm I\!I}^{+2}$ & 
0 & 
6 & 
8 & 
3 & 
1 & 
1 & 
$M_{24}$ 
\\  
 
90 & 
5 & 
24 & 
$S_4$& 
$2_{\rm I\!I}^{+4}8_{1}^{+1}3^{+2}$ & 
0 & 
2 & 
48 & 
1 & 
2 & 
2 & 
${\rm Mon}_a$ 
\\  
 
91 & 
5 & 
24 & 
$S_4$& 
$4_{5}^{+5}$ & 
0 & 
32 & 
64 & 
6 & 
1 & 
3 & 
$M_{24}$ 
\\  
 
92 & 
5 & 
24 & 
$S_4$& 
$2_{2}^{-4}8_{1}^{+1}3^{+1}$ & 
0 & 
1 & 
2 & 
1 & 
2 & 
1 & 
${\rm Mon}_b$ 
\\  
 
93 & 
5 & 
16 & 
$2^4$& 
$4_{4}^{+4}8_{1}^{+1}$ & 
0 & 
2 & 
256 & 
1 & 
1 & 
1 & 
${\rm Mon}_a$ 
\\  
 
94 & 
5 & 
16 & 
$[2^4]$\ (\#12)& 
$2_{4}^{+4}16_{1}^{+1}$ & 
0 & 
1 & 
1 & 
1 & 
2 & 
1 & 
${\rm Mon}_a$ 
\\  
 
95 & 
5 & 
16 & 
$[2^4]$\ (\#11)& 
$2_{2}^{+2}4_{1}^{+1}8_{2}^{+2}$ & 
0 & 
1 & 
16 & 
1 & 
2 & 
1 & 
${\rm Mon}_a$ 
\\  
 
96 & 
5 & 
12 & 
$[2^23]$\ (\#4)& 
$2_{\rm I\!I}^{+4}4_{1}^{+1}3^{-4}$ & 
0 & 
1 & 
1440 & 
1 & 
1 & 
- & 
${\rm Mon}_a$ 
\\  
 
97 & 
5 & 
12 & 
$[2^23]$\ (\#4)& 
$2_{4}^{+4}4_{1}^{+1}3^{-2}$ & 
0 & 
1 & 
4 & 
1 & 
3 & 
1 & 
${\rm Mon}_a$ 
\\  
 
98 & 
5 & 
12 & 
$[2^23]$\ (\#4)& 
$2_{2}^{-4}4_{1}^{+1}3^{-3}$ & 
0 & 
1 & 
48 & 
1 & 
3 & 
1 & 
${\rm Mon}_a$ 
\\  
 
\hline  
%%%%%%%%%%%%%%%%%%%%%%%%%%%%%%%%%%%%%%%%%%%%%%%%%%%%%%%%%%%%%%%%%%%%%% 

99 & 
4 & 
245760 & 
$2^8{:}M_{20}$& 
$2_{\rm I\!I}^{-2}4_{\rm I\!I}^{-2}$ & 
0 & 
1 & 
1 & 
1 & 
1 & 
- & 
${\rm Mon}_a$* 
\\  
 
100 & 
4 & 
30720 & 
$[2^9].A_5$& 
$2_{\rm I\!I}^{-4}5^{-1}$ & 
0 & 
1 & 
1 & 
1 & 
1 & 
- & 
${\rm Mon}_a$* 
\\  
 
101 & 
4 & 
29160 & 
$3^4{:}A_6$& 
$3^{+2}9^{+1}$ & 
1 & 
1 & 
1 & 
1 & 
1 & 
- & 
S* 
\\  
 
102 & 
4 & 
20160 & 
$L_3(4)$& 
$2_{\rm I\!I}^{-2}3^{-1}7^{-1}$ & 
2 & 
1 & 
1 & 
1 & 
2 & 
1 & 
$M_{23}$* 
\\  
 
103 & 
4 & 
12288 & 
$[2^{12}3]$& 
$2_{\rm I\!I}^{+2}4_{3}^{+1}8_{1}^{+1}$ & 
0 & 
1 & 
2 & 
1 & 
1 & 
- & 
${\rm Mon}_a$ 
\\  
 
104 & 
4 & 
9216 & 
$[2^{10}3^2]$& 
$2_{\rm I\!I}^{+4}3^{+2}$ & 
0 & 
1 & 
2 & 
1 & 
1 & 
- & 
${\rm Mon}_a$* 
\\  
 
105 & 
4 & 
6144 & 
$[2^{11}3]$& 
$2_{\rm I\!I}^{-2}4_{6}^{+2}3^{-1}$ & 
0 & 
1 & 
4 & 
1 & 
1 & 
- & 
${\rm Mon}_a$ 
\\  
 
106 & 
4 & 
5760 & 
$2^4{:}A_6$& 
$4_{5}^{-1}8_{1}^{+1}3^{+1}$ & 
2 & 
1 & 
1 & 
1 & 
2 & 
1 & 
$M_{23}$* 
\\  
 
107 & 
4 & 
4096 & 
$2^{1+8}{:}2^3$& 
$4_{4}^{+4}$ & 
0 & 
1 & 
8 & 
1 & 
1 & 
- & 
${\rm Mon}_a$ 
\\  
 
108 & 
4 & 
2520 & 
$A_7$& 
$3^{+1}5^{+1}7^{+1}$ & 
3 & 
1 & 
1 & 
1 & 
2 & 
1 & 
$M_{23}$* 
\\  
 
109 & 
4 & 
1944 & 
$3^{1+4}{:}2.2^2$& 
$2_{2}^{+2}3^{+3}$ & 
1 & 
1 & 
1 & 
1 & 
1 & 
- & 
S 
\\  
 
110 & 
4 & 
1920 & 
$2^4{:}S_5$& 
$4_{3}^{-1}8_{1}^{+1}5^{-1}$ & 
2 & 
1 & 
2 & 
1 & 
3 & 
1 & 
$M_{23}$ 
\\  
 
111 & 
4 & 
1344 & 
$2^3{:}L_2(7)$& 
$4_{2}^{+2}7^{+1}$ & 
2 & 
1 & 
1 & 
1 & 
3 & 
1 & 
$M_{23}$* 
\\  
 
112 & 
4 & 
1152 & 
$Q(3^2{:}2)$& 
$8_{6}^{-2}3^{-1}$ & 
2 & 
1 & 
2 & 
1 & 
2 & 
1 & 
$M_{23}$* 
\\  
 
113 & 
4 & 
1152 & 
$[2^73^2]$& 
$2_{\rm I\!I}^{-2}4_{6}^{+2}3^{-1}$ & 
0 & 
2 & 
4 & 
1 & 
1 & 
1 & 
$M_{24}$* 
\\  
 
114 & 
4 & 
972 & 
$[2^23^5]$\ (\#812)& 
$2_{\rm I\!I}^{-2}3^{+4}$ & 
0 & 
1 & 
6 & 
1 & 
1 & 
- & 
S 
\\  
 
115 & 
4 & 
768 & 
$[2^83]$\ & 
$2_{\rm I\!I}^{-2}8_{4}^{-2}$ & 
0 & 
6 & 
8 & 
3 & 
1 & 
1 & 
$M_{24}$ 
\\  
 
116 & 
4 & 
768 & 
$(2 \times 4^2){:}S_4$ & 
$2_{6}^{+2}8_{6}^{+2}$ & 
0 & 
1 & 
4 & 
1 & 
1 & 
- & 
${\rm Mon}_b$ 
\\  
 
117 & 
4 & 
768 & 
$[2^83]$\ & 
$4_{4}^{+4}$ & 
0 & 
12 & 
8 & 
6 & 
1 & 
1 & 
$M_{24}$ 
\\  
 
118 & 
4 & 
720 & 
$S_6$& 
$2_{\rm I\!I}^{-2}3^{+2}5^{+1}$ & 
2 & 
3 & 
3 & 
2 & 
2 & 
1 & 
$M_{23}$ 
\\  
 
119 & 
4 & 
720 & 
$M_{10}$& 
$2_{5}^{+1}4_{1}^{+1}3^{-1}5^{+1}$ & 
2 & 
1 & 
1 & 
1 & 
3 & 
2 & 
$M_{23}$ 
\\  
 
120 & 
4 & 
660 & 
$L_2(11)$& 
$11^{+2}$ & 
2 & 
1 & 
1 & 
1 & 
3 & 
2 & 
$M_{23}$* 
\\  
 
121 & 
4 & 
576 & 
$2^4{:}(S_3 \times S_3)$& 
$4_{7}^{+1}8_{1}^{+1}3^{+2}$ & 
2 & 
1 & 
4 & 
1 & 
3 & 
1 & 
$M_{23}$ 
\\  
 
122 & 
4 & 
500 & 
$5^{1+2}{:}4$& 
$5^{+3}$ & 
1 & 
1 & 
1 & 
1 & 
1 & 
- & 
S* 
\\  
 
123 & 
4 & 
384 & 
$[2^73]$\ (\#20097)& 
$2_{2}^{-2}4_{\rm I\!I}^{+2}3^{+1}$ & 
0 & 
1 & 
2 & 
1 & 
1 & 
- & 
${\rm Mon}_b$ 
\\  
 
124 & 
4 & 
384 & 
$[2^73]$\ (\#18134)& 
$2_{5}^{+3}16_{7}^{+1}$ & 
0 & 
1 & 
1 & 
1 & 
1 & 
- & 
${\rm Mon}_b$ 
\\  
 
125 & 
4 & 
384 & 
$[2^73]$\ (\#17948)& 
$2_{\rm I\!I}^{-2}4_{5}^{+1}8_{1}^{+1}3^{-1}$ & 
0 & 
3 & 
12 & 
2 & 
1 & 
1 & 
$M_{24}$ 
\\  
 
126 & 
4 & 
384 & 
$ [2^73]$\ (\#20089)& 
$2_{\rm I\!I}^{-2}4_{2}^{+2}3^{+1}$ & 
0 & 
6 & 
4 & 
3 & 
1 & 
1 & 
$M_{24}$ 
\\  
 
127 & 
4 & 
360 & 
$A_6$& 
$2_{\rm I\!I}^{+4}3^{+2}$ & 
0 & 
2 & 
2 & 
1 & 
1 & 
1 & 
${\rm Mon}_b$* 
\\  
 
128 & 
4 & 
360 & 
$(3\times A_5){:}2$& 
$3^{-2}5^{-2}$ & 
2 & 
1 & 
2 & 
1 & 
3 & 
2 & 
$M_{23}$* 
\\  
 
129 & 
4 & 
336 & 
$2 \times L_2(7)$& 
$2_{\rm I\!I}^{+2}7^{+2}$ & 
2 & 
1 & 
2 & 
1 & 
3 & 
2 & 
$M_{23}$ 
\\  
 
130 & 
4 & 
256 & 
$[2^8]$\ (\#53387)& 
$4_{3}^{+3}8_{1}^{+1}$ & 
0 & 
1 & 
16 & 
1 & 
1 & 
- & 
${\rm Mon}_a$ 
\\  
 
131 & 
4 & 
192 & 
$[2^63]$\ (\#1494)& 
$2_{6}^{+2}8_{6}^{+2}$ & 
0 & 
2 & 
4 & 
1 & 
1 & 
1 & 
$M_{24}$* 
\\  
 
132 & 
4 & 
192 & 
$[2^63]$\ (\#1493)& 
$4_{4}^{+4}$ & 
0 & 
12 & 
8 & 
3 & 
1 & 
2 & 
$M_{24}$* 
\\  
 
133 & 
4 & 
144 & 
$[2^43^2]$\ (\#187)& 
$2_{4}^{+4}9^{-1}$ & 
0 & 
1 & 
1 & 
1 & 
2 & 
1 & 
${\rm Mon}_a$ 
\\  
 
134 & 
4 & 
144 & 
$3^2{:}QD_{16}$& 
$2_{1}^{+1}4_{1}^{+1}3^{-1}9^{-1}$ & 
2 & 
1 & 
2 & 
1 & 
4 & 
2 & 
$M_{23}$ 
\\  
 
135 & 
4 & 
144 & 
$[2^43^2]$\ (\#183)& 
$2_{\rm I\!I}^{-2}4_{\rm I\!I}^{-2}3^{-2}$ & 
0 & 
16 & 
32 & 
6 & 
1 & 
2 & 
$M_{24}$ 
\\  
 
136 & 
4 & 
120 & 
$S_5$& 
$2_{\rm I\!I}^{+4}5^{-2}$ & 
0 & 
2 & 
12 & 
1 & 
1 & 
1 & 
${\rm Mon}_b$ 
\\  
 
137 & 
4 & 
120 & 
$S_5$& 
$2_{6}^{+4}3^{+1}5^{-1}$ & 
0 & 
1 & 
2 & 
1 & 
2 & 
1 & 
${\rm Mon}_b$ 
\\  
 
138 & 
4 & 
108 & 
$[2^23^3]$\ (\#17)& 
$3^{+2}9^{+2}$ & 
0 & 
18 & 
54 & 
4 & 
1 & 
3 & 
$M_{24}$* 
\\  
 
139 & 
4 & 
96 & 
$[2^53]$\ (\#195)& 
$2_{\rm I\!I}^{+2}4_{\rm I\!I}^{+2}3^{+2}$ & 
0 & 
2 & 
16 & 
1 & 
1 & 
1 & 
${\rm Mon}_a$ 
\\  
 
140 & 
4 & 
72 & 
$[2^33^2]$\ (\#46)& 
$2_{\rm I\!I}^{+4}3^{+2}9^{-1}$ & 
0 & 
1 & 
72 & 
1 & 
1 & 
- & 
${\rm Mon}_a$ 
\\  
 
141 & 
4 & 
72 & 
$[2^33^2]$\ (\#40)& 
$2_{2}^{-4}3^{+1}9^{-1}$ & 
0 & 
1 & 
6 & 
1 & 
2 & 
1 & 
${\rm Mon}_b$ 
\\  
 
142 & 
4 & 
72 & 
$[2^43^2]$\ (\#40)& 
$2_{\rm I\!I}^{+4}3^{+2}9^{-1}$ & 
0 & 
2 & 
72 & 
1 & 
1 & 
1 & 
${\rm Mon}_b$ 
\\  
 
143 & 
4 & 
64 & 
$[2^6]$\ (\#257)& 
$2_{3}^{+1}4_{1}^{+1}8_{\rm I\!I}^{+2}$ & 
0 & 
1 & 
8 & 
1 & 
1 & 
- & 
${\rm Mon}_b$ 
\\  
 
144 & 
4 & 
60 & 
$A_5$& 
$2_{\rm I\!I}^{-2}3^{+4}$ & 
0 & 
6 & 
6 & 
2 & 
1 & 
2 & 
$M_{24}$* 
\\  
 
145 & 
4 & 
48 & 
$[2^43]$\ (\#48)& 
$2_{\rm I\!I}^{+2}4_{7}^{+1}8_{1}^{+1}3^{+2}$ & 
0 & 
1 & 
48 & 
1 & 
2 & 
1 & 
${\rm Mon}_a$ 
\\  
 
146 & 
4 & 
48 & 
$[2^43]$\ (\#48)& 
$2_{\rm I\!I}^{-2}4_{1}^{+1}8_{1}^{+1}3^{+1}$ & 
0 & 
1 & 
4 & 
1 & 
2 & 
1 & 
${\rm Mon}_a$ 
\\  
 
147 & 
4 & 
48 & 
$[2^43]$\ (\#29)& 
$4_{\rm I\!I}^{-2}8_{\rm I\!I}^{-2}$ & 
0 & 
32 & 
64 & 
6 & 
1 & 
3 & 
${\rm Mon}_b$* 
\\  
 
148 & 
4 & 
48 & 
$[2^43]$\ (\#48)& 
$2_{0}^{-2}4_{1}^{+1}8_{1}^{+1}3^{+1}$ & 
0 & 
1 & 
4 & 
1 & 
2 & 
1 & 
${\rm Mon}_a$ 
\\  
 
149 & 
4 & 
40 & 
$[2^35]$\ (\#12)& 
$2_{4}^{+4}5^{+2}$ & 
0 & 
3 & 
6 & 
2 & 
2 & 
3 & 
$M_{24}$ 
\\  
 
150 & 
4 & 
36 & 
$[2^23^2]$\ (\#13)& 
$2_{4}^{+4}3^{+4}$ & 
0 & 
3 & 
72 & 
2 & 
1 & 
1 & 
${\rm Mon}_a$ 
\\  
 
151 & 
4 & 
36 & 
$[2^23^2]$\ (\#10)& 
$2_{6}^{-4}3^{+3}$ & 
0 & 
3 & 
6 & 
2 & 
1 & 
1 & 
${\rm Mon}_b$ 
\\  
 
152 & 
4 & 
32 & 
$[2^5]$\ (\#40)& 
$2_{2}^{+2}4_{1}^{+1}16_{1}^{+1}$ & 
0 & 
1 & 
2 & 
1 & 
2 & 
1 & 
${\rm Mon}_a$ 
\\  
 
153 & 
4 & 
32 & 
$[2^5]$\ (\#46)& 
$4_{2}^{+2}8_{2}^{+2}$ & 
0 & 
2 & 
64 & 
1 & 
1 & 
1 & 
${\rm Mon}_a$ 
\\  
 
154 & 
4 & 
24 & 
$[2^33]$\ (\#14)& 
$2_{2}^{+2}4_{2}^{+2}3^{-2}$ & 
0 & 
1 & 
8 & 
1 & 
2 & 
1 & 
${\rm Mon}_a$ 
\\  
 
155 & 
4 & 
24 & 
$[2^33]$\ (\#8)& 
$2_{\rm I\!I}^{-2}4_{\rm I\!I}^{-2}3^{+4}$ & 
0 & 
1 & 
1152 & 
1 & 
1 & 
- & 
${\rm Mon}_a$ 
\\  
 
156 & 
4 & 
24 & 
$2^2.S_3$& 
$2_{\rm I\!I}^{-2}4_{\rm I\!I}^{-2}3^{+4}$ & 
0 & 
6 & 
1152 & 
2 & 
1 & 
2 & 
$M_{24}$* 
\\  
 
157 & 
4 & 
24 & 
$[2^33]$\ (\#6)& 
$2_{6}^{-2}4_{\rm I\!I}^{+2}3^{+3}$ & 
0 & 
1 & 
48 & 
1 & 
1 & 
- & 
${\rm Mon}_b$ 
\\  
 
158 & 
4 & 
24 & 
$[2^33]$\ (\#14)& 
$2_{\rm I\!I}^{-2}4_{2}^{+2}3^{-3}$ & 
0 & 
1 & 
96 & 
1 & 
1 & 
- & 
${\rm Mon}_a$ 
\\  
 
159 & 
4 & 
20 & 
$[2^25]$\ (\#4)& 
$2_{\rm I\!I}^{-4}5^{-3}$ & 
0 & 
1 & 
120 & 
1 & 
1 & 
- & 
${\rm Mon}_a$ 
\\  
 
160 & 
4 & 
16 & 
$[2^4]$\ (\#3)& 
$4_{2}^{+2}8_{2}^{+2}$ & 
0 & 
4 & 
64 & 
1 & 
1 & 
3 & 
$M_{24}$ 
\\  
 
161 & 
4 & 
12 & 
$[2^23]$\ (\#4)& 
$2_{4}^{+4}3^{+4}$ & 
0 & 
18 & 
72 & 
4 & 
1 & 
5 & 
$M_{24}$ 
\\  
 
\hline  
%%%%%%%%%%%%%%%%%%%%%%%%%%%%%%%%%%%%%%%%%%%%%%%%%%%%%%%%%%%%%%%%%%%%%% 

162 & 
3 & 
10321920 & 
$2^9.L_3(4)$& 
$2_{\rm I\!I}^{-2}8_{3}^{-1}$ & 
0 & 
1 & 
1 & 
1 & 
1 & 
- & 
${\rm Mon}_a$* 
\\  
 
163 & 
3 & 
3265920 & 
$U_4(3)$& 
$4_{7}^{+1}3^{+2}$ & 
1 & 
1 & 
1 & 
1 & 
1 & 
- & 
S* 
\\  
 
164 & 
3 & 
491520 & 
$[2^{12}].S_5$& 
$4_{3}^{+3}$ & 
0 & 
1 & 
2 & 
1 & 
1 & 
- & 
${\rm Mon}_a$ 
\\  
 
165 & 
3 & 
443520 & 
$M_{22}$& 
$4_{5}^{-1}11^{+1}$ & 
2 & 
1 & 
1 & 
1 & 
2 & 
1 & 
$M_{23}$* 
\\  
 
166 & 
3 & 
184320 & 
$[2^9].A_6$& 
$2_{\rm I\!I}^{-2}4_{1}^{+1}3^{+1}$ & 
0 & 
1 & 
1 & 
1 & 
1 & 
- & 
${\rm Mon}_a$* 
\\  
 
167 & 
3 & 
126000 & 
$U_3(5)$& 
$2_{7}^{+1}5^{-2}$ & 
1 & 
1 & 
1 & 
1 & 
1 & 
- & 
S* 
\\  
 
168 & 
3 & 
61440 & 
$[2^9].S_5$& 
$2_{\rm I\!I}^{-2}4_{7}^{+1}5^{-1}$ & 
0 & 
1 & 
2 & 
1 & 
1 & 
- & 
${\rm Mon}_a$ 
\\  
 
169 & 
3 & 
58320 & 
$3^4.A_6.2$& 
$2_{1}^{+1}3^{-1}9^{+1}$ & 
1 & 
1 & 
1 & 
1 & 
1 & 
- & 
S 
\\  
 
170 & 
3 & 
40320 & 
$L_3(4).2$& 
$4_{3}^{-1}3^{-1}7^{-1}$ & 
2 & 
1 & 
2 & 
1 & 
2 & 
1 & 
$M_{23}$ 
\\  
 
171 & 
3 & 
40320 & 
$L_3(4).2$& 
$2_{5}^{+3}7^{-1}$ & 
0 & 
1 & 
1 & 
1 & 
1 & 
- & 
- 
\\  
 
172 & 
3 & 
40320 & 
$2^4.A_7$& 
$8_{1}^{+1}7^{+1}$ & 
2 & 
1 & 
1 & 
1 & 
2 & 
1 & 
$M_{23}$* 
\\  
 
173 & 
3 & 
36864 & 
$[2^{12}3^2]$& 
$2_{\rm I\!I}^{+2}8_{5}^{-1}3^{-1}$ & 
0 & 
1 & 
2 & 
1 & 
1 & 
- & 
${\rm Mon}_a$* 
\\  
 
174 & 
3 & 
24576 & 
$[2^{13}3]$& 
$4_{2}^{+2}8_{1}^{+1}$ & 
0 & 
1 & 
4 & 
1 & 
1 & 
- & 
${\rm Mon}_a$ 
\\  
 
175 & 
3 & 
20160 & 
$A_8$& 
$4_{1}^{+1}3^{+1}5^{+1}$ & 
2 & 
1 & 
1 & 
1 & 
2 & 
1 & 
$M_{23}$* 
\\  
 
176 & 
3 & 
18432 & 
$[2^{11}3^2]$& 
$2_{\rm I\!I}^{+2}4_{7}^{+1}3^{+2}$ & 
0 & 
1 & 
4 & 
1 & 
1 & 
- & 
${\rm Mon}_a$ 
\\  
 
177 & 
3 & 
12288 & 
$[2^{12}3]$& 
$4_{1}^{-3}3^{-1}$ & 
0 & 
1 & 
8 & 
1 & 
1 & 
- & 
${\rm Mon}_a$ 
\\  
 
178 & 
3 & 
11520 & 
$2^4.S_6$& 
$2_{2}^{-2}8_{7}^{+1}3^{+1}$ & 
0 & 
1 & 
2 & 
1 & 
1 & 
- & 
${\rm Mon}_b$ 
\\  
 
179 & 
3 & 
11520 & 
$2^4.S_6$& 
$2_{\rm I\!I}^{-2}8_{1}^{+1}3^{+1}$ & 
0 & 
6 & 
4 & 
3 & 
1 & 
1 & 
$M_{24}$ 
\\  
 
180 & 
3 & 
10752 & 
$[2^6].L_2(7)$& 
$2_{2}^{-2}16_{5}^{-1}$ & 
0 & 
1 & 
1 & 
1 & 
1 & 
- & 
${\rm Mon}_b$* 
\\  
 
181 & 
3 & 
10752 & 
$[2^6].L_2(7)$& 
$4_{3}^{+3}$ & 
0 & 
2 & 
2 & 
1 & 
1 & 
1 & 
$M_{24}$* 
\\  
 
182 & 
3 & 
7920 & 
$M_{11}$& 
$2_{7}^{+1}3^{-1}11^{-1}$ & 
2 & 
1 & 
1 & 
1 & 
3 & 
2 & 
$M_{23}$* 
\\  
 
183 & 
3 & 
5760 & 
$[2^43].S_5$& 
$8_{1}^{+1}3^{-1}5^{-1}$ & 
2 & 
1 & 
2 & 
1 & 
3 & 
1 & 
$M_{23}$* 
\\  
 
184 & 
3 & 
4608 & 
$[2^93^2]$& 
$2_{3}^{-1}8_{\rm I\!I}^{-2}$ & 
0 & 
1 & 
2 & 
1 & 
1 & 
- & 
${\rm Mon}_b$* 
\\  
 
185 & 
3 & 
4608 & 
$[2^93^2]$& 
$4_{2}^{+2}8_{1}^{+1}$ & 
0 & 
2 & 
4 & 
1 & 
1 & 
1 & 
$M_{24}$* 
\\  
 
186 & 
3 & 
3888 & 
$[2^43^5]$& 
$4_{1}^{+1}3^{+3}$ & 
0 & 
1 & 
2 & 
1 & 
1 & 
- & 
S 
\\  
 
187 & 
3 & 
3888 & 
$[2^43^5]$& 
$2_{3}^{+3}3^{-2}$ & 
0 & 
1 & 
1 & 
1 & 
1 & 
- & 
S 
\\  
 
188 & 
3 & 
3840 & 
$[2^5].S_5$& 
$2_{\rm I\!I}^{-2}8_{7}^{+1}5^{-1}$ & 
0 & 
3 & 
6 & 
2 & 
1 & 
1 & 
$M_{24}$ 
\\  
 
189 & 
3 & 
2688 & 
$[2^4].L_2(7)$& 
$2_{\rm I\!I}^{+2}4_{1}^{+1}7^{+1}$ & 
0 & 
3 & 
3 & 
2 & 
1 & 
1 & 
$M_{24}$ 
\\  
 
190 & 
3 & 
2304 & 
$[2^83^2]$& 
$4_{1}^{-3}3^{-1}$ & 
0 & 
6 & 
8 & 
3 & 
1 & 
1 & 
$M_{24}$ 
\\  
 
191 & 
3 & 
1944 & 
$[2^33^5]$\ (\#3536)& 
$2_{1}^{-3}3^{-3}$ & 
0 & 
1 & 
6 & 
1 & 
1 & 
- & 
S 
\\  
 
192 & 
3 & 
1536 & 
$[2^93]$& 
$4_{1}^{+1}8_{2}^{+2}$ & 
0 & 
1 & 
8 & 
1 & 
1 & 
- & 
${\rm Mon}_a$ 
\\  
 
193 & 
3 & 
1440 & 
$2.A_6.2$& 
$2_{2}^{-2}4_{1}^{+1}5^{+1}$ & 
0 & 
1 & 
1 & 
1 & 
2 & 
1 & 
${\rm Mon}_a$ 
\\  
 
194 & 
3 & 
1440 & 
$A_6.2^2$& 
$2_{5}^{+3}3^{-1}5^{+1}$ & 
0 & 
3 & 
3 & 
2 & 
1 & 
1 & 
$M_{24}$ 
\\  
 
195 & 
3 & 
1152 & 
$[2^73^2]$& 
$2_{\rm I\!I}^{+2}8_{7}^{+1}3^{+2}$ & 
0 & 
1 & 
8 & 
1 & 
1 & 
- & 
${\rm Mon}_a$ 
\\  
 
196 & 
3 & 
768 & 
$[2^83]$\ & 
$4_{5}^{-3}3^{+1}$ & 
0 & 
1 & 
4 & 
1 & 
1 & 
- & 
${\rm Mon}_a$ 
\\  
 
197 & 
3 & 
768 & 
$[2^83]$\ & 
$2_{3}^{+1}4_{1}^{+1}16_{7}^{+1}$ & 
0 & 
1 & 
2 & 
1 & 
1 & 
- & 
${\rm Mon}_b$ 
\\  
 
198 & 
3 & 
720 & 
$A_6.2$& 
$2_{\rm I\!I}^{+2}4_{7}^{+1}3^{+2}$ & 
0 & 
2 & 
4 & 
1 & 
1 & 
1 & 
${\rm Mon}_b$ 
\\  
 
199 & 
3 & 
720 & 
$A_6.2$& 
$2_{\rm I\!I}^{+2}4_{1}^{+1}3^{+1}5^{+1}$ & 
0 & 
1 & 
6 & 
1 & 
2 & 
1 & 
${\rm Mon}_b$ 
\\  
 
200 & 
3 & 
720 & 
$A_6.2$& 
$2_{6}^{+2}4_{1}^{+1}3^{+2}$ & 
0 & 
1 & 
2 & 
1 & 
1 & 
- & 
${\rm Mon}_b$ 
\\  
 
201 & 
3 & 
432 & 
$[2^43^3]$\ (\#734)& 
$2_{7}^{+1}3^{+2}9^{-1}$ & 
0 & 
6 & 
6 & 
2 & 
1 & 
2 & 
$M_{24}$* 
\\  
 
202 & 
3 & 
384 & 
$[2^73]$\ (\#5602)& 
$4_{1}^{+1}8_{2}^{+2}$ & 
0 & 
2 & 
8 & 
1 & 
1 & 
1 & 
$M_{24}$ 
\\  
 
203 & 
3 & 
288 & 
$[2^53^2]$\ (\#1027)& 
$2_{2}^{+2}4_{1}^{+1}9^{-1}$ & 
0 & 
1 & 
2 & 
1 & 
2 & 
1 & 
${\rm Mon}_a$ 
\\  
 
204 & 
3 & 
240 & 
$2.S_5$& 
$2_{\rm I\!I}^{+2}4_{7}^{+1}5^{-2}$ & 
0 & 
1 & 
12 & 
1 & 
1 & 
- & 
${\rm Mon}_a$ 
\\  
 
205 & 
3 & 
240 & 
$2.S_5$& 
$2_{4}^{+2}4_{1}^{+1}3^{+1}5^{-1}$ & 
0 & 
1 & 
4 & 
1 & 
2 & 
1 & 
${\rm Mon}_a$ 
\\  
 
206 & 
3 & 
192 & 
$[2^63]$\ (\#1472)& 
$4_{7}^{+3}3^{+2}$ & 
0 & 
1 & 
16 & 
1 & 
1 & 
- & 
${\rm Mon}_a$ 
\\  
 
207 & 
3 & 
168 & 
$[2^33.7]$\ (\#43)& 
$2_{3}^{-1}8_{\rm I\!I}^{-2}$ & 
0 & 
2 & 
2 & 
1 & 
1 & 
1 & 
${\rm Mon}_b$* 
\\  
 
208 & 
3 & 
144 & 
$[2^43^2]$\ (\#183)& 
$2_{2}^{+2}8_{1}^{+1}3^{-2}$ & 
0 & 
1 & 
4 & 
1 & 
1 & 
- & 
${\rm Mon}_b$ 
\\  
 
209 & 
3 & 
144 & 
$[2^43^2]$\ (\#186)& 
$2_{\rm I\!I}^{-2}4_{1}^{+1}3^{+1}9^{-1}$ & 
0 & 
1 & 
12 & 
1 & 
1 & 
- & 
${\rm Mon}_a$ 
\\  
 
210 & 
3 & 
144 & 
$[2^43^2]$\ (\#189)& 
$2_{\rm I\!I}^{-2}8_{1}^{+1}3^{-3}$ & 
0 & 
1 & 
48 & 
1 & 
1 & 
- & 
${\rm Mon}_a$ 
\\  
 
211 & 
3 & 
120 & 
$S_5$& 
$2_{1}^{-3}3^{-3}$ & 
0 & 
6 & 
6 & 
2 & 
1 & 
2 & 
$M_{24}$ 
\\  
 
212 & 
3 & 
96 & 
$[2^53]$\ (\#189)& 
$2_{4}^{-2}16_{1}^{+1}3^{-1}$ & 
0 & 
1 & 
2 & 
1 & 
2 & 
1 & 
${\rm Mon}_a$ 
\\  
 
213 & 
3 & 
96 & 
$[2^53]$\ (\#226)& 
$4_{4}^{-2}8_{1}^{+1}3^{+1}$ & 
0 & 
2 & 
16 & 
1 & 
1 & 
1 & 
${\rm Mon}_a$ 
\\  
 
214 & 
3 & 
72 & 
$[2^33^2]$\ (\#46)& 
$2_{4}^{-2}4_{1}^{+1}3^{+3}$ & 
0 & 
3 & 
12 & 
2 & 
1 & 
1 & 
${\rm Mon}_a$ 
\\  
 
215 & 
3 & 
64 & 
$ [2^6]$\ (\#131)& 
$4_{2}^{+2}16_{1}^{+1}$ & 
0 & 
2 & 
8 & 
1 & 
1 & 
1 & 
${\rm Mon}_a$ 
\\  
 
216 & 
3 & 
64 & 
$[2^6]$\ (\#73)& 
$8_{3}^{+3}$ & 
0 & 
2 & 
16 & 
1 & 
1 & 
1 & 
${\rm Mon}_a$ 
\\  
 
217 & 
3 & 
48 & 
$[2^43]$\ (\#48)& 
$2_{\rm I\!I}^{-2}8_{1}^{+1}3^{-3}$ & 
0 & 
6 & 
48 & 
2 & 
1 & 
2 & 
$M_{24}$ 
\\  
 
218 & 
3 & 
48 & 
$[2^43]$\ (\#51)& 
$4_{3}^{+3}3^{-2}$ & 
0 & 
1 & 
16 & 
1 & 
1 & 
- & 
${\rm Mon}_a$ 
\\  
 
219 & 
3 & 
48 & 
$[2^43]$\ (\#38)& 
$4_{1}^{-3}3^{+3}$ & 
0 & 
1 & 
96 & 
1 & 
1 & 
- & 
${\rm Mon}_a$ 
\\  
 
220 & 
3 & 
42 & 
$[2.3.7]$\ (\#1)& 
$2_{3}^{+3}7^{-2}$ & 
0 & 
1 & 
6 & 
1 & 
1 & 
- & 
- 
\\  
 
221 & 
3 & 
24 & 
$[2^33]$\ (\#3)& 
$8_{3}^{+3}$ & 
0 & 
4 & 
16 & 
1 & 
1 & 
3 & 
${\rm Mon}_b$* 
\\  
 
\hline  
%%%%%%%%%%%%%%%%%%%%%%%%%%%%%%%%%%%%%%%%%%%%%%%%%%%%%%%%%%%%%%%%%%%%%% 
222 & 
2 & 
9196830720 & 
$U_6(2)$& 
$2_{\rm I\!I}^{-2}3^{+1}$ & 
0 & 
1 & 
1 & 
1 & 
1 & 
- & 
S* 
\\  
 
223 & 
2 & 
898128000 & 
${\rm McL}$& 
$3^{-1}5^{-1}$ & 
1 & 
1 & 
1 & 
1 & 
1 & 
- & 
S* 
\\  
 
224 & 
2 & 
454164480 & 
$2^{10}.M_{22}$& 
$4_{2}^{+2}$ & 
0 & 
1 & 
1 & 
1 & 
1 & 
- & 
${\rm Mon}_a$* 
\\  
 
225 & 
2 & 
44352000 & 
${\rm HS}$& 
$2_{2}^{-2}5^{+1}$ & 
0 & 
1 & 
1 & 
1 & 
1 & 
- & 
S* 
\\  
 
226 & 
2 & 
20643840 & 
$2^9.L_3(4).2$& 
$4_{1}^{+1}8_{1}^{+1}$ & 
0 & 
1 & 
2 & 
1 & 
1 & 
- & 
${\rm Mon}_a$ 
\\  
 
227 & 
2 & 
10200960 & 
$M_{23}$& 
$23^{+1}$ & 
1 & 
1 & 
1 & 
1 & 
2 & 
1 & 
$M_{23}$* 
\\  
 
228 & 
2 & 
6531840 & 
$U_4(3).2$& 
$2_{6}^{+2}3^{+2}$ & 
0 & 
1 & 
2 & 
1 & 
1 & 
- & 
S 
\\  
 
229 & 
2 & 
6531840 & 
$U_4(3).2$& 
$2_{3}^{-1}4_{1}^{+1}3^{-1}$ & 
0 & 
1 & 
1 & 
1 & 
1 & 
- & 
S 
\\  
 
230 & 
2 & 
1924560 & 
$3^5.M_{11}$& 
$3^{-1}9^{-1}$ & 
0 & 
1 & 
1 & 
1 & 
1 & 
- & 
S* 
\\  
 
231 & 
2 & 
1474560 & 
$[2^{12}3].S_5$& 
$4_{\rm I\!I}^{-2}3^{-1}$ & 
0 & 
1 & 
2 & 
1 & 
1 & 
- & 
${\rm Mon}_a$* 
\\  
 
232 & 
2 & 
1290240 & 
$[2^9].A_7$& 
$2_{\rm I\!I}^{+2}7^{+1}$ & 
0 & 
1 & 
1 & 
1 & 
1 & 
- & 
${\rm Mon}_a$* 
\\  
 
233 & 
2 & 
887040 & 
$M_{22}.2$& 
$2_{\rm I\!I}^{-2}11^{+1}$ & 
0 & 
3 & 
3 & 
2 & 
1 & 
1 & 
$M_{24}$ 
\\  
 
234 & 
2 & 
368640 & 
$[2^9].A_6.2$& 
$4_{4}^{-2}3^{+1}$ & 
0 & 
1 & 
2 & 
1 & 
1 & 
- & 
${\rm Mon}_b$ 
\\  
 
235 & 
2 & 
322560 & 
$2^4.A_8$& 
$2_{3}^{-1}16_{3}^{-1}$ & 
0 & 
1 & 
1 & 
1 & 
1 & 
- & 
${\rm Mon}_b$* 
\\  
 
236 & 
2 & 
322560 & 
$2^4.A_8$& 
$4_{1}^{+1}8_{1}^{+1}$ & 
0 & 
2 & 
2 & 
1 & 
1 & 
1 & 
$M_{24}$* 
\\  
 
237 & 
2 & 
184320 & 
$[2^93].S_5$& 
$2_{\rm I\!I}^{+2}3^{-1}5^{-1}$ & 
0 & 
1 & 
2 & 
1 & 
1 & 
- & 
${\rm Mon}_a$* 
\\  
 
238 & 
2 & 
147456 & 
$[2^93^2]$& 
$8_{2}^{+2}$ & 
0 & 
1 & 
2 & 
1 & 
1 & 
- & 
${\rm Mon}_a$* 
\\  
 
239 & 
2 & 
122880 & 
$[2^{10}].S_5$& 
$4_{2}^{-2}5^{-1}$ & 
0 & 
1 & 
4 & 
1 & 
1 & 
- & 
${\rm Mon}_a$ 
\\  
 
240 & 
2 & 
120960 & 
$L_3(4).[2.3]$& 
$3^{+2}7^{-1}$ & 
0 & 
2 & 
4 & 
1 & 
1 & 
1 & 
$M_{24}$* 
\\  
 
241 & 
2 & 
116640 & 
$3^4.2.A_6.2$& 
$2_{2}^{+2}9^{+1}$ & 
0 & 
1 & 
1 & 
1 & 
1 & 
- & 
S 
\\  
 
242 & 
2 & 
95040 & 
$M_{12}$& 
$2_{6}^{+2}3^{+2}$ & 
0 & 
2 & 
2 & 
1 & 
1 & 
1 & 
$M_{24}$* 
\\  
 
243 & 
2 & 
80640 & 
$L_3(4).2^2$& 
$2_{3}^{+1}4_{1}^{+1}7^{-1}$ & 
0 & 
1 & 
2 & 
1 & 
1 & 
- & 
- 
\\  
 
244 & 
2 & 
73728 & 
$[2^{13}3^2]$& 
$4_{7}^{-1}8_{1}^{+1}3^{-1}$ & 
0 & 
1 & 
4 & 
1 & 
1 & 
- & 
${\rm Mon}_a$ 
\\  
 
245 & 
2 & 
58320 & 
$3^4.A_6.2$& 
$2_{\rm I\!I}^{-2}3^{+1}9^{+1}$ & 
0 & 
1 & 
6 & 
1 & 
1 & 
- & 
S 
\\  
 
246 & 
2 & 
40320 & 
$2.A_8$& 
$2_{\rm I\!I}^{+2}3^{+1}5^{+1}$ & 
0 & 
1 & 
2 & 
1 & 
1 & 
- & 
${\rm Mon}_a$ 
\\  
 
247 & 
2 & 
36864 & 
$[2^{12}3^2]$& 
$4_{6}^{+2}3^{+2}$ & 
0 & 
1 & 
8 & 
1 & 
1 & 
- & 
${\rm Mon}_a$ 
\\  
 
248 & 
2 & 
23040 & 
$[2^5].A_6.2$& 
$4_{3}^{-1}8_{1}^{+1}3^{+1}$ & 
0 & 
1 & 
4 & 
1 & 
1 & 
- & 
${\rm Mon}_a$ 
\\  
 
249 & 
2 & 
23040 & 
$[2^63].S_5$& 
$4_{2}^{-2}5^{-1}$ & 
0 & 
2 & 
4 & 
1 & 
1 & 
1 & 
$M_{24}$* 
\\  
 
250 & 
2 & 
21504 & 
$[2^7].L_2(7)$& 
$4_{1}^{+1}16_{1}^{+1}$ & 
0 & 
1 & 
2 & 
1 & 
1 & 
- & 
${\rm Mon}_b$ 
\\  
 
251 & 
2 & 
15840 & 
$2\times M_{11}$& 
$2_{4}^{-2}11^{-1}$ & 
0 & 
1 & 
1 & 
1 & 
2 & 
1 & 
${\rm Mon}_a$ 
\\  
 
252 & 
2 & 
7920 & 
$M_{11}$& 
$4_{\rm I\!I}^{-2}3^{-1}$ & 
0 & 
2 & 
2 & 
1 & 
1 & 
1 & 
${\rm Mon}_b$* 
\\  
 
253 & 
2 & 
7776 & 
$[2^53^5]$& 
$2_{1}^{+1}4_{1}^{+1}3^{-2}$ & 
0 & 
1 & 
2 & 
1 & 
1 & 
- & 
S 
\\  
 
254 & 
2 & 
6912 & 
$[2^83^3]$& 
$4_{6}^{+2}3^{+2}$ & 
0 & 
2 & 
8 & 
1 & 
1 & 
1 & 
$M_{24}$* 
\\  
 
255 & 
2 & 
5040 & 
$S_7$& 
$2_{\rm I\!I}^{-2}5^{+1}7^{+1}$ & 
0 & 
1 & 
6 & 
1 & 
1 & 
- & 
- 
\\  
 
256 & 
2 & 
5040 & 
$S_7$& 
$2_{6}^{-2}3^{+1}7^{+1}$ & 
0 & 
1 & 
2 & 
1 & 
1 & 
- & 
- 
\\  
 
257 & 
2 & 
2880 & 
$2.A_6.2^2$& 
$4_{6}^{-2}5^{+1}$ & 
0 & 
2 & 
4 & 
1 & 
1 & 
1 & 
${\rm Mon}_b$ 
\\  
 
258 & 
2 & 
2688 & 
$[2^4].L_2(7)$& 
$8_{2}^{+2}$ & 
0 & 
2 & 
2 & 
1 & 
1 & 
1 & 
${\rm Mon}_b$* 
\\  
 
259 & 
2 & 
2304 & 
$[2^83^2]$& 
$2_{5}^{-1}16_{7}^{+1}3^{-1}$ & 
0 & 
1 & 
2 & 
1 & 
1 & 
- & 
${\rm Mon}_b$ 
\\  
 
260 & 
2 & 
2160 & 
$3.A_6.2$& 
$2_{\rm I\!I}^{-2}3^{+1}9^{+1}$ & 
0 & 
6 & 
6 & 
2 & 
1 & 
2 & 
$M_{24}$* 
\\  
 
261 & 
2 & 
1440 & 
$2.A_6.2$& 
$4_{6}^{+2}3^{+2}$ & 
0 & 
2 & 
8 & 
1 & 
1 & 
1 & 
${\rm Mon}_b$ 
\\  
 
262 & 
2 & 
1000 & 
$5^{1+2}.[2^3]$\ (\#91)& 
$2_{2}^{+2}5^{+2}$ & 
0 & 
1 & 
2 & 
1 & 
1 & 
- & 
S 
\\  
 
263 & 
2 & 
864 & 
$[2^53^3]$\ (\#4661)& 
$2_{4}^{-2}3^{-1}9^{-1}$ & 
0 & 
3 & 
3 & 
2 & 
1 & 
1 & 
${\rm Mon}_a$ 
\\  
 
264 & 
2 & 
720 & 
$[2.3].S_5$\ (\#767)& 
$2_{\rm I\!I}^{-2}3^{-1}5^{-2}$ & 
0 & 
1 & 
12 & 
1 & 
1 & 
- & 
- 
\\  
 
265 & 
2 & 
720 & 
$[2.3].S_5$\ (\#767)& 
$2_{6}^{-2}3^{-2}5^{-1}$ & 
0 & 
1 & 
4 & 
1 & 
1 & 
- & 
- 
\\  
 
266 & 
2 & 
576 & 
$[2^63^2]$\ (\#8299)& 
$4_{2}^{+2}9^{-1}$ & 
0 & 
1 & 
4 & 
1 & 
1 & 
- & 
${\rm Mon}_a$ 
\\  
 
267 & 
2 & 
432 & 
$ [2^43^3]$\ (\#523)& 
$4_{\rm I\!I}^{-2}3^{-1}9^{-1}$ & 
0 & 
1 & 
12 & 
1 & 
1 & 
- & 
${\rm Mon}_a$* 
\\  
 
268 & 
2 & 
432 & 
$[2^43^3]$\ (\#734)& 
$4_{\rm I\!I}^{-2}3^{-1}9^{-1}$ & 
0 & 
2 & 
12 & 
1 & 
1 & 
1 & 
${\rm Mon}_b$* 
\\  
 
269 & 
2 & 
384 & 
$[2^73]$\ (\#18127)& 
$8_{1}^{+1}16_{1}^{+1}$ & 
0 & 
2 & 
4 & 
1 & 
1 & 
1 & 
${\rm Mon}_a$* 
\\  
 
270 & 
2 & 
288 & 
$[2^53^2]$\ (\#1028)& 
$4_{1}^{+1}8_{1}^{+1}3^{-2}$ & 
0 & 
1 & 
8 & 
1 & 
1 & 
- & 
${\rm Mon}_a$ 
\\  
 
271 & 
2 & 
240 & 
$2.S_5$& 
$4_{6}^{+2}3^{+2}$ & 
0 & 
4 & 
8 & 
1 & 
1 & 
3 & 
$M_{24}$* 
\\  
 
272 & 
2 & 
80 & 
$[2^45]$\ (\#34)& 
$4_{2}^{+2}5^{+2}$ & 
0 & 
1 & 
8 & 
1 & 
1 & 
- & 
${\rm Mon}_a$ 
\\  
 
\hline  
%%%%%%%%%%%%%%%%%%%%%%%%%%%%%%%%%%%%%%%%%%%%%%%%%%%%%%%%%%%%%%%%%%%%%% 
273 & 
1 & 
$|{\rm Co}_2|$ & 
${\rm Co}_2$ & 
$4_{1}^{+1}$ & 
0 & 
1 & 
1 & 
1 & 
1 & 
- & 
S* 
\\  
 
274 & 
1 & 
$|{\rm Co}_3|$ & 
${\rm Co}_3$& 
$2_{3}^{-1}3^{-1}$ & 
0 & 
1 & 
1 & 
1 & 
1 & 
- & 
S* 
\\  
 
275 & 
1 & 
20891566080 & 
$ 2^{11}.M_{23}$& 
$8_{1}^{+1}$ & 
0 & 
1 & 
1 & 
1 & 
1 & 
- & 
${\rm Mon}_a*$ 
\\  
 
276 & 
1 & 
18393661440 & 
$ U_6(2).2$& 
$4_{3}^{-1}3^{+1}$ & 
0 & 
1 & 
2 & 
1 & 
1 & 
- & 
S 
\\  
 
277 & 
1 & 
1796256000 & 
${\rm McL}.2$& 
$2_{5}^{-1}5^{-1}$ & 
0 & 
1 & 
1 & 
1 & 
1 & 
- & 
S 
\\  
 
278 & 
1 & 
244823040 & 
$M_{24}$& 
$4_{3}^{-1}3^{+1}$ & 
0 & 
2 & 
2 & 
1 & 
1 & 
1 & 
$M_{24}$* 
\\  
 
279 & 
1 & 
88704000 & 
${\rm HS}.2$& 
$4_{5}^{-1}5^{+1}$ & 
0 & 
1 & 
2 & 
1 & 
1 & 
- & 
S 
\\  
 
280 & 
1 & 
61931520 & 
$2^9.L_3(4).3.2$& 
$8_{3}^{-1}3^{-1}$ & 
0 & 
1 & 
2 & 
1 & 
1 & 
- & 
${\rm Mon}_a$* 
\\  
 
281 & 
1 & 
10321920 & 
$2^9.A_8$& 
$16_{1}^{+1}$ & 
0 & 
1 & 
1 & 
1 & 
1 & 
- & 
${\rm Mon}_b$* 
\\  
 
282 & 
1 & 
3849120 & 
$3^5.(2 \times M_{11})$& 
$2_{1}^{+1}9^{-1}$ & 
0 & 
1 & 
1 & 
1 & 
1 & 
- & 
S 
\\  
 
283 & 
1 & 
2580480 & 
$2^9.S_7$& 
$4_{7}^{+1}7^{+1}$ & 
0 & 
1 & 
2 & 
1 & 
1 & 
- & 
- 
\\  
 
284 & 
1 & 
737280 & 
$2^{1+4+6}.3.S_5$& 
$8_{5}^{-1}5^{-1}$ & 
0 & 
1 & 
2 & 
1 & 
1 & 
- & 
${\rm Mon}_a$* 
\\  
 
285 & 
1 & 
368640 & 
$2^{10}.(3\times A_5).2$& 
$4_{7}^{+1}3^{-1}5^{-1}$ & 
0 & 
1 & 
4 & 
1 & 
1 & 
- & 
- 
\\  
 
286 & 
1 & 
241920 & 
$L_3(4).3.2.2$& 
$2_{5}^{-1}3^{-1}7^{-1}$ & 
0 & 
1 & 
2 & 
1 & 
1 & 
- & 
- 
\\  
 
287 & 
1 & 
233280 & 
$3^4.2.A_6.2.2$& 
$4_{1}^{+1}9^{+1}$ & 
0 & 
1 & 
2 & 
1 & 
1 & 
- & 
S 
\\  
 
288 & 
1 & 
190080 & 
$ 2.M_{12}$& 
$8_{3}^{-1}3^{-1}$ & 
0 & 
2 & 
2 & 
1 & 
1 & 
1 & 
${\rm Mon}_b$* 
\\  
 
289 & 
1 & 
3456 & 
$(2^3 \times 3^2).Q_8.S_3$& 
$8_{1}^{+1}9^{-1}$ & 
0 & 
1 & 
2 & 
1 & 
1 & 
- & 
${\rm Mon}_a$* 
\\  
 
\hline  
%%%%%%%%%%%%%%%%%%%%%%%%%%%%%%%%%%%%%%%%%%%%%%%%%%%%%%%%%%%%%%%%%%%%%% 

290 & 
1 & 
$|{\rm Co}_0|$ & 
${\rm Co}_0$ & 
$1$ & 
0 & 
1 & 
1 & 
1 & 
1 & 
- & 
S* 
\\  

\hline  
%%%%%%%%%%%%%%%%%%%%%%%%%%%%%%%%%%%%%%%%%%%%%%%%%%%%%%%%%%%%%%%%%%%%%% 

\end{longtable} }

%%%%%%%%%%%%%%%%%%%%%%%%%%%%%%%%%%%%%%%%%%%%%%%%%%%%%%%%%
%\newpage
%
%\bigskip
%
%\sloppy
%
%{\footnotesize
%
%\renewcommand{\arraystretch}{0.71}
%\arraycolsep0.85mm
%
%\begin{table}[h]
%\caption{Gram matrices of the fixed-point lattices.}\label{gram}
%\phantom{xxx}
%\end{table}
%\input{data/table2.erg}}
%
%%%%%%%%%%%%%%%%%%%%%%%%%%%%%%%%%%%%%%%%%%%%%%%%%%%%%%%%%
%
%\bigskip
%
%{\footnotesize
%
%\renewcommand{\arraystretch}{0.71}
%\arraycolsep0.85mm
%
%\begin{table}[h]
%\caption{Lattice of fixed-point lattice stabilizers. For each orbit, the list of orbits determined by
%subgroups of the stabilizer which are maximal among the subgroups with a larger fixed-point lattice is given.
%A $*$ after the list indicates that not all subgroups have been determined and there may orbits be missing.}\label{stablat}
%\phantom{xxx}
%\end{table}
%
%\noindent\input{data/table3.erg}}
%
%}
%
%%%%%%%%%%%%%%%%%%%%%%%%%%%%%%%%%%%%%%%%%%%%%%%%%%%%%%%%%%%%%%%%%%%%%%%%%%%%%%%%%%%%%%%%%%%%%%%%%%%%%%%%%%%%%%%%%

\bigskip

We collect now several observations regarding the tables.\ In some cases
these may be read-off directly from the tables, while others can be obtained
by simple arguments or easy calculations.\ In any case we omit details.

\medskip

%%%%%%%%%%%%

\paragraph{The isometry type of the lattices $\Lambda^G$ and $\Lambda_G$.}
The isometry class of the coinvariant lattice $\Lambda_G$ determines uniquely the orbit of $\Lambda^G$.\
However, isometric $\Lambda^G$ may belong to different orbits.\
In the following table we itemize the \emph{isometric orbits\/} (i.e., orbits of isometric fixed-point lattices) 
which contain more than one orbit of fixed-point lattices.

\smallskip

\noindent$\begin{array}{r|l}
\mbox{Rank} & \mbox{Sets of isometric lattices $\Lambda^G$} \\ \hline
6 &    \{ 34, 36 \},\ \{ 40, 49 \},\ \{ 41, 48, 56 \},\  \{ 57, 64 \} \\[2mm]
5 &    \{ 67, 71 \},\ \{ 78, 86, 91 \} \\[2mm]
4 &    \{ 104, 127 \},\ \{ 105, 113 \},\ \{ 107, 117, 132 \},\  \{ 114, 144 \},\  \{ 116, 131 \},\  \{ 140, 142 \},\  \\
  & \qquad  \{ 150, 161 \},\ \{ 153, 160 \} \\[2mm]
3 &    \{ 155, 156 \},\ \{ 164, 181 \},\  \{ 174, 185 \},\ \{ 176, 198 \},\ \{ 177, 190 \},\ \{ 184, 207 \},\\\
  &    \qquad  \{ 191, 211 \},\   \{ 192, 202 \},\   \{ 210, 217 \},\  \{ 216, 221 \} \\[2mm]
2 &    \{ 226, 236 \},\   \{ 228, 242 \},\ \{ 231, 252 \},\ \{ 238, 258 \},\ \{ 239, 249 \},\ \{ 245, 260 \},\  \\
  & \qquad \{ 247, 254, 261, 271 \},\ \{ 267, 268 \} \\[2mm]
1 &   \{ 276, 278 \},\ \{ 280, 288 \} 
\end{array}$

The lattices $\Lambda^G$ and $\Lambda_G$ are isometric to each other in all three rank~$12$ cases.

\paragraph{The genus of $\Lambda^G$ and $\Lambda_G$.}

The genera of $\Lambda_G$ and $\Lambda^G$ determine each other.\ 
Two orbits of fixed-point lattices $\Lambda^G$ define the same genus if, and only if, they are isometric.

The isometry classes of lattices in the genus of $\Lambda^G$ have the following property: if the class
belongs to a fixed-point lattice then the minimal norm is at least~$4$; for all other classes, 
the root sublattice has maximal rank.\ The root lattice of $\Lambda^G$ itself is zero exactly
for orbits no.\ $1$, $2$, $4$, $7$, $18$, $20$, $39$, $52$, $53$, $82$, $108$, $120$, $128$, $129$, $227$, $243$, $251$. 
These lattices were investigated (without explicit classification) 
in~\cite{Bo-Leech}.\ Most of them are fixed-point lattices of conjugacy classes in $M_{23}$. 

As for the isometry classes of lattices in the genus of $\Lambda_G$, if the class
belongs to $\Lambda_G$ then the minimal norm is~$4$. For all other classes it seems that 
the minimal norm is $2$ although the root lattice does \emph{not\/} always has maximal rank.
However, we checked this only in a small number of cases.

%%%%%%%%%%%%

\paragraph{The entry $\alpha$.}
For an even lattice $L$ we define $\alpha(L)= {\rm rk}\, L - {\rm rk}\,A_L$.\
Clearly \hbox{$\alpha(L)\ge 0$.}

1.\ We have $\alpha(\Lambda^G)\ge 2$ if, and only if, $G\subsetneq M_{23}$, i.e., $G=\widetilde G$ 
is a proper subgroup of the stabilizer of lattice no.\ 227.

2.\ We have $\alpha(\Lambda^G)\geq 1$ if, and only if, $G\subseteq\mcl$ or $G\subseteq M_{23}$, 
i.e., $G$ is contained in the stabilizer of either lattice no.\ 223 or lattice no.\ 227.

%%%%%%%%%%%%

\paragraph{Niemeier lattices.}
Let $N$ be a Niemeier lattice in the sense that it is one of the $24$ lattices in the genus of $\Lambda$.\
Its isometry group is a split extension $O(N) =W(N){:}G$,
where $W(N)$ is generated by reflections in hyperplanes orthogonal to the roots of $N$.\
The coinvariant lattice $N_G$, which is always a lattice without norm~$2$ vectors, can be embedded into $\Lambda$
in such a way that $G\cong O_0(N_G) \cong O_0(\Lambda_G)$ (cf.~\cite{Nikulin}, Remark~1.14.7, Prop.~1.14.8 and~\cite{Conway3}).
The following table lists the no.\ of the corresponding entry of $\Lambda_G$ in Table~\ref{maintable}.

\smallskip

\noindent {\footnotesize
$\begin{array}{l|*{12}{r}}
{\rm Lattice} & D_{24} & E_8^3 & D_{16}E_8 & A_{24} & D_{12}^2 & D_{10}E_7^2 & A_{17}E_7 & A_{15}D_9 & D_8^3 & A_{12}^2 & A_{11}D_7E_6 & E_6^4 \\ \hline
{\rm No.} & 1 & 22  & 1 & 5 & 5 & 2 & 2  & 2  & 22 & 64 & 2 & 147
\end{array}$ }

\noindent {\footnotesize
$\begin{array}{l|*{12}{r}}
{\rm Lattice} & D_6^4  &  A_9^2 D_6 & A_8^3 & A_7^2D_5^2 & A_6^4 & D_4^6 & A_5^4D_4 & A_4^6 & A_3^8 & A_2^{12} & A_1^{24} & \Lambda  \\ \hline
{\rm No.}     &  91    &     9      & 161   &     21     & 221   & 260   &   87     & 271   & 258   & 288      &  278     & 290                                                                 
\end{array}$ }

\medskip

Conversely, to obtain all embeddings of a given $\Lambda_G$ from Table~\ref{maintable} into Niemeier lattices with roots, 
we determined all isometry classes of lattices $K$ in the genus of $\Lambda^G$ and all equivalence classes 
of extensions $K\oplus \Lambda_G$ to an unimodular lattice $N$.\
There is always a \emph{unique\/} lattice $K$ providing a \emph{unique extension\/} of $K \oplus \Lambda_G$ to the Leech lattice $\Lambda$.
Column $N$ of Table~\ref{maintable} lists the number of isomorphism classes of Niemeier lattices $N$ with roots obtained in this
way.\ If this number is positive, $G$ embeds into the group $O(N)/W(N)$ of the corresponding Niemeier lattices $N$.

%%%%%%%%%%%%

\paragraph{Conjugacy classes of $\Co$.}
There are $72$ conjugacy classes $[g]$ in $\Co$ such that $\Lambda^g\not= 0$, giving rise
to $58$ fixed-point lattices $\Lambda^{\langle g\rangle}$ considered in~\cite{HaLa}.\ Below we list 
these lattices, their rank, and the index of the image of $N_{\Co}(\langle g\rangle)$ in $O(\Lambda^{\langle g\rangle})$.

\smallskip

\noindent{\footnotesize 
$\begin{array}{l|*{18}{r}}
\mbox{order} & 1 & 2 & 2 & 2 & 3 & 3 & 3 & 4 & 4 & 4 & 4 & 4 & 4 & 5 & 5 & 6 & 6 & 6  \\ \hline
\mbox{rank}    & 24& 8 & 16 & 12 & 12 & 6 & 8 & 8 & 6 & 10 & 4 & 8 & 6 & 8 & 4 & 6 & 6 & 6 \\ \hline
\mbox{no.}     & 1 & 14 & 2 & 5 & 4 & 35 & 22 & 14 & 41 & 9 & 99 & 21 & 64 & 20 & 122 & 35 & 62 & 33  \\ \hline
\mbox{index}   & 1 & 2 & 1 & 5040 & 1 & 1 & 1920 & 240 & 1 & 1 & 2 & 36 & 6 & 1 & 1 & 1 & 1 & 1  \\ \hline
\end{array}$ }

\smallskip

\noindent{\footnotesize 
$\begin{array}{l|*{18}{r}}
\mbox{order} & 6 & 6 & 6 & 6 & 6 & 6 & 7 & 8 & 8 & 8 & 8 & 8 & 8 & 9 & 9 & 10 & 10 & 10  \\ \hline
\mbox{rank}    & 4 &  8 & 4 & 2 & 6 & 4 & 6 & 4 & 4 & 2 & 6 & 4 & 4 & 2 & 4 & 4 & 4 & 4 \\ \hline
\mbox{no.}     & 104 &  18 & 114 & 222 & 63 & 161 & 52 & 99 & 107 & 224 & 55 & 143 & 147 & 230 & 101 & 100 & 122 & 159 \\ \hline
\mbox{index}   & 2 & 1 & 2 & 1 & 1 & 4 & 1 & 6 & 2 & 1 & 1 & 1 & 6 & 1 & 2 & 1 & 1 & 1 \\ \hline
\end{array}$ }

\smallskip

\noindent{\footnotesize 
$\begin{array}{l|*{18}{r}}
\mbox{order} & 10 & 11 & 12 & 12 & 12 & 12 & 12 & 12 & 12 & 12 & 12 & 12 & 14 & 14 & 15 & 15 & 16 & 16 \\ \hline
\mbox{rank}    &  4 & 4 & 2 & 4 & 2 & 2 & 2 & 4 & 4 & 4 & 4 & 2 & 4 & 2 & 4 & 2 & 2 & 2 \\ \hline
\mbox{no.}     &  149 & 120 & 222 & 104 & 228 & 222 & 231 & 109 & 123 & 157 & 135 & 271 & 129 & 232 & 128 & 223 & 224 & 226 \\ \hline
\mbox{index}   &  2 & 1 & 1 & 2 & 1 & 1 & 1 & 1 & 1 & 1 & 6 & 1 & 1 & 1 & 1 & 1 & 1 & 1 \\ \hline
\end{array}$ }

\smallskip

\noindent{\footnotesize 
$\begin{array}{l|*{18}{r}}
\mbox{order} & 18 & 18 & 18 & 20 & 20 & 20 & 21 & 22 & 22 & 23 & 23 & 24 & 24 & 24 & 28 & 30 & 30 & 30\\ \hline
\mbox{rank}    &  2 & 2 & 2 & 2 & 2 & 2 & 2 & 2 & 2 & 2 & 2 & 2 & 2 & 2 & 2 & 2 & 2 & 2 \\ \hline
\mbox{no.}     & 230 & 222 & 245 & 262 & 257 & 225 & 240 & 251 & 251 & 227 & 227 & 229 & 234 & 253 & 232 & 237 & 223 & 246 \\ \hline
\mbox{index}   & 1 & 1 & 1 & 1 & 1 & 1 & 1 & 1 & 1 & 1 & 1 & 1 & 1 & 1 & 1 & 1 & 1 & 1 \\ \hline
\end{array}$ }

%%%%%%%%%%%%

\paragraph{$\S$-Lattices.}\label{slattices}
Each of the twelve $\S$-lattices~\cite{Curtis} $S$ arises as a fixed-point lattice in $\Lambda$.\
The type of $S$, denoted by $2^a3^b$, records the numbers $a$, $b$ of pairs of short 
representatives $\pm v$ for $S/2S$ of norm~$4$, $6$ respectively.\ 
For an $\S$-lattice we always have $1+a+b=2^{\rk(S)}$ and
$S$ is characterized up to isometry by its type.\ 
The $\S$-lattices are identified in the following table.

\smallskip

\noindent {\footnotesize
$\begin{array}{l|*{12}{r}}
\S\mbox{-Lattice} & 2^03^0 & 2^13^0 & 2^03^1 & 2^33^0 & 2^23^1 & 2^13^2 & 2^03^3 & 2^53^2 & 2^33^4 & 2^93^6 & 2^53^{10} & 2^{27}3^{36} \\ \hline
{\rm rank\ }         &    0   &    1   &   1    &   2    &   2    &   2    &    2   &  3     &   3    &    4   &     4     &   6      \\ \hline 
{\rm no.}            &  290   &   273  & 274    &  222   &  223   &  225   &  230   &  163   &  167   &  101   &  122      &  35
\end{array}$ }

\smallskip

The stabilizer $H$ of some $\S$-lattices can be extended to a stabilizer $G$ with a lower 
dimensional non-trivial fixed-point lattice $S'$ such that $H=\widetilde{O^2(G)}$.
The following $18$ orbits arise:

\smallskip

\noindent {\footnotesize
$\begin{array}{l|r|r|r|r|rr|rrrr|r}
\S \mbox{-Lattice} & 2^33^0 & 2^23^1 & 2^13^2 & 2^03^3 & \multicolumn{2}{c|}{2^53^2} &  \multicolumn{4}{c|}{2^93^6}     & 2^53^{10} \\ \hline
|G/H| &               2     &   2    &   2     &  2     & 2 &  2    &  2 & 2 & 2^2 & 2^3       &     2     \\  \hline
{\rm rank\ } S'  &  1     &   1    &   1     &   1    & 2 &  2    &  3 & 2 &  2 &   1        &     2     \\ \hline
{\rm no.}        & 276    & 277    & 279     &   282  & 228 &229  &  169 & 245 & 241 & 287   &    262                      
\end{array}$ }

\medskip

\noindent {\footnotesize
$\begin{array}{l|*{8}{r}}
\S\mbox{-Lattice} &   \multicolumn{7}{c}{2^{27}3^{36}} \\ \hline
|G/H|  &          2  & 2^2 & 2   & 2^3  & 2^3  &  2^2 & 2^4   \\  \hline
{\rm rank\ } S'  &   5 & 4   & 4   &  3   & 3    & 3    & 2      \\ \hline
{\rm No.}        &   68 & 109 & 114 &  186 &  187 &  191 & 253                   
\end{array}$ }

%%%%%%%%%%%%

\paragraph{Groups related to $\mon$.}
Let $G$ be the full stabilizer of a lattice such that $\Lambda^{\widetilde{O^2(G)}}$ is not an $\S$-lattice.\
Using  inclusions $M_{23}\subseteq M_{24}\subseteq \mon \subseteq \Co$,
we define the following five \emph{types\/} of $G$:
\vspace{-2mm}
\begin{itemize}
\itemsep0pt
\item[$M_{23}$:] $G$ is contained in $M_{23}$ ($61$ cases);
\item[$M_{24}$:] $G$ is contained in $M_{24}$ but not in $M_{23}$  ($128-61=67$ cases);
\item[${\rm Mon}_a$:] $G$ is contained in $\mon$ but not in $M_{24}$ and $G=T{:}H$ with 
$H\subseteq M_{24}$ and $T=G\cap 2^{12}$ ($212-128 = 84$ cases);
\item[${\rm Mon}_b$:] $G$ is contained in $\mon$ but not of type ${\rm Mon}_a$  ($250-212=38$ cases);
\item[$-$:] $G$ is not contained in $\mon$  ($10$ cases).
\end{itemize}
The type of each $G$ is listed in the last column of Table~\ref{maintable}.

If $H\subseteq M_{23}$ then $\widetilde{H}\subseteq M_{23}$.\ For $H$ is contained
in $2^{11}.M_{23}$ and $M_{24}$, which are both stabilizers of rank~$1$ lattices, whence 
(with an obvious notation) $\widetilde{H}\subseteq 2^{11}.M_{23}\cap M_{24}= M_{23}$.\ Similarly,
$H\subseteq M_{24}$ implies $\widetilde{H}\subseteq M_{24}$.\ If $H\subseteq \mon$ but is contained in neither $2^{11}.M_{23}$ nor $M_{24}$, then
$\widetilde{H}$ is  generally \emph{not\/} contained in $\mon$.

%%%%%%%%%%%%

\paragraph{Spherical Designs.}
The even integral lattices of minimal norm~$4$ for which the minimal vectors form 
spherical $6$-designs have been classified by Martinet~\cite{Martinet}. All of them can be obtained 
from~$\Lambda$.\ In the nomenclature of Table~\ref{maintable}
they are as follows:\ $2\Z$ ($\Lambda^G$ no.\ 273), $E_8(2)$ ($\Lambda^G$ no.\ 14 or $\Lambda_G$ no.\ 2),
the Barnes-Wall lattice of rank $16$ ($\Lambda^G$ no.\ 2 or $\Lambda_G$ no.\ 14), 
$\Lambda_{23}$ ($\Lambda_G$ no.\ 273), and $\Lambda$ itself.

A lattice whose minimal vectors and those of its dual form spherical $4$-designs is called 
\emph{dual strongly perfect.\/}\
Using the Molien series of their full automorphism groups, the following additional
lattices can be shown to be dual strongly perfect, cf.~\cite{Venkov}:
$A_2$ ($\Lambda^G$ for no.\ 222),
$D_4$ ($\Lambda^G$ for no.\ 99),
$E_6$ ($\Lambda^G$ for no.\ 33, 35),
one lattice of rank $10$ ($\Lambda^G$ for no.\ 7),
Coxeter-Todd lattice $K_{12}$ ($\Lambda^G\cong \Lambda_G$ for no.\ 4),
one lattice of rank~$18$ ($\Lambda_G$ for no.\ 35),
two lattices of rank~$22$ ($\Lambda_G$ for no.\ 222 and no.\ 223),
one lattice of rank~$23$ ($\Lambda_G$ for no.\ 274).\

In addition, further lattices $\Lambda^G$ and $\Lambda_G$ are rescaled versions 
of the above listed lattices.

%%%%%%%%%%%%%%%%%%%%%%%%%%%%%%%%%%%%%%%%%%%%%%%%%%%%%%%%%%%%%%%%%%%%%%%%%%%%%%%%%%%%%%%%%

\end{document}